\documentclass[10pt,reqno,final]{amsart} 
\usepackage{epsfig,amssymb,amsmath,version}
\usepackage{amssymb,version,graphicx,fancybox,mathrsfs,mathtools}
\usepackage[notcite,notref]{showkeys}  
\usepackage{url,hyperref}
\usepackage{subfigure}
\usepackage{color}
\usepackage{stmaryrd}
\usepackage{multirow}
\usepackage{booktabs,siunitx}
\usepackage{booktabs}
\usepackage{multicol}
\usepackage{float}   
\usepackage{algorithm}
\usepackage{algorithmic}
\usepackage{setspace}
\usepackage{array}
\usepackage{booktabs} 
\usepackage{multirow}
\usepackage{diagbox}
\usepackage{epstopdf}
\usepackage{cite}
\usepackage{makecell}

\catcode`\@=11 \theoremstyle{plain}
\@addtoreset{equation}{section}   

\@addtoreset{figure}{section}
\renewcommand\thefigure{\thesection.\@arabic\c@figure}
\newtheorem{thm}{\bf Theorem}

\newtheorem{cor}{\bf Corollary}

\newtheorem{lmm}{\bf Lemma}

\theoremstyle{remark}
\newtheorem{rem}{\bf Remark}[section]

\definecolor{ligreen}{rgb}{0.0, 0.3, 0.0}

\definecolor{darkblue}{rgb}{0.0, 0.0, 0.55}

\definecolor{anti-flashwhite}{rgb}{0.55, 0.57, 0.68}

\newcommand{\bs}[1]{\boldsymbol{#1}}

\graphicspath{{./Figures/} } 

\pdfoutput=1

\begin{document}

\baselineskip 14pt

\bibliographystyle{plain}
\title[An Efficient Algorithm for Multiple Solutions] {An
Efficient Spectral  Trust-Region Deflation Method for Multiple Solutions}
\author[
	L. Li,\,    L. Wang\,  $\&$\,  H. Li
	]{
		\;\; Lin Li${}^{\dag}$,   \;\;  Li-Lian Wang${}^{\ddag}$ \;\; and\;\; Huiyuan Li${}^{\S}$
		}
	\thanks{${}^{\dag}$School of Mathematics and Physics and Hunan Key Laboratory for Mathematical Modeling and Scientific Computing, University of South China, Hengyang, 421001, China. This work of this author is partially supported by the Science Foundations of Hunan Province (Nos: 2020JJ5464, 20C1595). Email: lilinmath@usc.edu.cn (L. Li).\\
			\indent ${}^{\ddag}$Corresponding author.  Division of Mathematical Sciences, School of Physical and Mathematical Sciences, Nanyang Technological University,
		637371, Singapore. The research of this author is partially supported by Singapore MOE AcRF Tier 1 Grant: RG15/21.  Email: lilian@ntu.edu.sg (L. Wang).\\
		\indent ${}^{\S}${State Key Laboratory of Computer Science/Laboratory of Parallel Computing, Institute of Software, Chinese Academy of Sciences, Beijing 100190, China}. Email: huiyuan@iscas.ac.cn (H. Li).\\[2pt]
\indent The first author would like to thank the hospitality of Nanyang Technological University for the visit to finalize this work.}


\begin{abstract} It is quite common that a nonlinear partial differential equation (PDE) admits multiple distinct solutions and each solution may carry a unique physical meaning. One typical approach for finding multiple solutions is to  use the Newton method with different initial guesses that ideally fall into the basins of attraction confining the solutions.   In this paper,  we propose a fast and accurate numerical method for multiple solutions comprised of three ingredients: (i) a well-designed spectral-Galerkin discretization of the underlying PDE leading to a nonlinear  algebraic system  (NLAS) with multiple solutions; (ii) an effective deflation technique to eliminate a known (founded) solution from the other unknown solutions  leading to deflated NLAS; and (iii) a viable nonlinear least-squares and trust-region (LSTR) method  for solving the NLAS and the deflated NLAS to find the multiple solutions sequentially one  by one.
We demonstrate through ample examples of differential equations and comparison with relevant existing approaches that  the spectral LSTR-Deflation  method has the merits: (i) it is quite flexible in choosing initial values, even starting from the same initial guess for finding all multiple solutions; (ii) it guarantees high-order accuracy;  and  (iii) it is quite fast to locate multiple distinct solutions and explore new solutions  which are not reported in literature.
\end{abstract}
\maketitle

\section{Introduction}

Many nonlinear differential equations arisen from mathematical modeling in  physics, mechanics, biology, energy and engineering, admit multiple solutions
and each  may carry its own physics and unique dynamics (see, e.g., \cite{1995Introduction, 2002A, Tadmor2012A}). Most problems, if not all,  do not have explicit solutions, so the development of efficient and accurate  numerical methods for finding  multiple solutions becomes a topic of longstanding interest.

In this paper, we intend to identify multiple solutions of the nonlinear ordinary differential equations
\begin{equation}\label{New1.1}
F(x,u,u',\cdots,u^{(n)}) = 0,\quad x\in \Omega,
\end{equation}
and the second-order nonlinear elliptic problems
\begin{equation}\label{New1.2}
\Delta u + F(x, y, u) = 0,   \quad  (x, y) \in \Omega,
\end{equation}
supplemented with some boundary conditions, where $\Omega$ is a bounded domain in ${\mathbb R}^{d}, d=1,2,$
and the nonlinear functional  $F$ has a regularity to be specified later.  We refer to  \cite{Lions1982On, 1999Exact, 2003Multiple, 2007On}
for the study of the existence and multiplicity of the solutions to such problems.

Before we introduce our approach,  we feel compelled to elaborate on the relevant background and motivations. Choi and Mckenna \cite{choi1993mountain}   proposed the mountain pass algorithm (MPA) for (\ref{New1.2}). As commented in
Xie, Chen and Yu \cite{xie2005improved}, the MPA is feasible for  finding two solutions of mountain pass type with Morse index 1 or 0. However,
according to Ding, Costa and Chen \cite{ding1999high},  the MPA  may  fail to locate the sign-changing solutions, so they proposed
a high linking algorithm (HLA) to compute such solutions with some inspirations from \cite{wang1991superlinear}.
Later, Li and Zhou \cite{li2001minimax} developed  the minimax algorithm (MNA) for multiple solutions for which we refer to \cite{2005A, 2007Numerical, 2008Numerical, 2005Instability, 2011A} for more recent advancements along this line.
It is noteworthy that these methods are more or less based on the variational structures of the underlying nonlinear differential equations.
The second category of the existing methods is to discretize the differential equation by  a numerical method (e.g., finite difference, finite element  or spectral method) and then search for multiple solutions of the resulting nonlinear algebraic system (NLAS).
 In this regard, a search-extension method  (SEM) was proposed in  \cite{chen2004search} for the semilinear PDEs with some improvements in \cite{xie2006improved, xie2015augmented}.
 Roughly speaking, the SEM searches for the starter from solving the companion linear problem (e.g., by the eigenfunction expansion method)
 and then extend to solve the NLAS, e.g., by Newton iteration with the starter as initial guess.
  Noticeably, a homotopy continuation method was introduced  by Allgower et al. \cite{allgower2006solution, allgower2009application}
 for finding multiple solutions of the NLAS obtained from finite difference discretization of type \eqref{New1.2}.
 In general, in order to solve the NLAS:  $\bs F(\bs x)=\bs 0,$ we consider   ${\bs H}_\tau({\bs x}) = (1-\tau) {\bs F}({\bs x}) +
 \gamma\tau {\bs G}({\bs x})=\bs 0,$  where ${\bs H}_\tau({\bs x})$ is known as  the homotopy function,
  $\tau\in [0,1]$ is  the homotopy tracking number, and $\gamma$ is a randomly chosen complex number. Here, the starting system ${\bs H}_1({\bs x})=\bs G(\bs x) =\bs 0$ should be constructed which is expected to be easier to solve and closely related to the original NLAS. Then  the homotopy systems
  ${\bs H}_{\tau_j}({\bs x})=\bs 0$ should be solved for various $\tau_j$ from $1$ to $0$ using a suitable iterative method, where $\tau$ can be viewed as an extra dimension. This technique has inspired some interesting recent works, see, e.g., \cite{zhang2013eigenfunction,hao2014bootstrapping,2018Two}.

This work is mostly motivated by Farrell et al. \cite{farrell2015deflation}  where  a deflation technique and Newton iteration  were  integrated to solve the NLAS resulted from finite element discretization of the underlying PDE with multiple distinct  solutions. Given a solution  $\bs x_1$ of $\bs F(\bs x)=\bs 0,$
the essence of deflation is to eliminate this known solution $\bs x_1$ and then continue to search for a new solution by modifying the NLAS as $\bs G(\bs x; \bs x_1)=\bs 0,$ where
$\bs G$ is properly constructed so that  $\bs G(\bs x_1; \bs x_1)\not =\bs 0$. In a nutshell,  the algorithm in  \cite{farrell2015deflation}
 repeated this process and used the Newton iteration for finding $\bs x_1$ and solving the deflated systems.
In fact, the deflation technique  has its root in  polynomial root-finding problems  \cite{1966Rounding} and been extended  by Brown and Gearhart  \cite{1971Deflation} for solutions of  nonlinear algebraic systems.
As we know,  the  Newton iteration is sensitive to the initial guess and costly to compute the inverse of the Jacobian matrix that become even severer when the NLAS are deflated.   Indeed, the choice of good initial inputs becomes ad hoc and  problematic, and at times, this integrated algorithm fails to converge.


We aim to overcome the sensitivity of initial guesses and substantially  improve the efficiency and accuracy of the deflation technique.  
Our proposal consists of the following  ingredients.
\begin{itemize}
\item[(i)] We employ well-designed   Legendre spectral-Galerkin methods to discretise
\eqref{New1.1} and \eqref{New1.2} that can guarantee  a good approximability
of the resulted NLAS to the original differential equations with a relatively low computational cost.
\smallskip
\item[(ii)]  We make use of  the deflation technique  to sequentially  search for multiple solutions of the NLAS.
\smallskip
\item[(iii)] We introduce the nonlinear least-squares trust-region method to solve the NLAS and the deflated systems.
\end{itemize}
Accordingly, the proposed algorithm is dubbed as the spectral LSTR-Deflation method.  Compared with the  existing methods, the main differences and advantages of our algorithm reside in the following aspects.
\begin{itemize}
  \item Finding multiple solutions based on  the Newton method typically requires  the attempt of many different initial guesses that lie in different basins of attraction \cite{farrell2015deflation}.
   However, the use of the LSTR method has a good convergence where the choices of  initial inputs are much more relaxed and
     fairly flexible at times.  As the deflated systems turn out to be more and more complicated, the LSTR method becomes
     vitally important to ensure a good performance of  the whole algorithm.  This allows us to start even with the same initial guesses for multiple solutions.
\smallskip
  \item Compare with several existing methods, the proposed approach  is capable of finding new solutions more efficiently and more quickly (see the examples in Section \ref{sect3}).
\end{itemize}

\smallskip

The remainder of this paper is organized as follows. In Section \ref{sect2}, we describe the spectral LSTR-Deflation method for the model problem  (\ref{New1.2}).  In Section \ref{sect3}, we provide ample numerical experiments to demonstrate the efficiency  and accuracy of the algorithm. We then conclude the paper with some remarks.

\section{The Legendre spectral-Galerkin LSTR-Deflation method}\label{sect2}

In this section, we present the spectral LSTR-Deflation method for the model problem  (\ref{New1.2}). We start with its Legendre spectral-Galerkin discretisation   that leads to the nonlinear linear algebraic system to be solved. We then introduce the nonlinear least-squares trust-region method for the NLAS, followed by the exposition of the deflation technique. Finally, we summarise the whole algorithm.

\subsection{Legendre spectral-Galerkin method for a model problem}\label{subsect21}  As an illustrative example, we consider the nonlinear elliptic problem:
\begin{equation}\label{modelexample}
-\Delta u + F(u) = f(x,y),   \quad  (x, y)\in\Omega:=(-1,1)^2;\quad u|_{\partial \Omega}=0,
\end{equation}
where $F(z)$ is a smooth nonlinear function (e.g., a polynomial as in \cite{hao2018Two}), and   $f$ is continuous.

As the first step, we discretize \eqref{modelexample}  by using the Legendre spectral-Galerkin method    \cite{shen2011spectral}.
Let ${\mathbb P}_N$ be the set of the polynomials of degree at most $N,$ and ${\mathbb P}_N^0=\{\phi\in {\mathbb P}_N\,:\, \phi(\pm 1)=0\}.$
The spectral-Galerkin approximation is to find $u_N\in X_N^0=({\mathbb P}_N^0)^2$ such that
 \begin{equation}\label{New2.3B}
(\nabla u_N, \nabla v_N)+(I_N F(u_N), v_N) =  (I_{N}f, v_{N}), \quad  \forall v_{N}\in X_{N}^0,
\end{equation}
where $I_N$ the Legendre-Gauss-Lobatto tensorial interpolation operator with $N+1$ points in each coordinate direction. As in \cite{shen2011spectral}, we introduce the basis of ${\mathbb P}_N^0$ as follows
\begin{equation*}
  \phi_{k}(x) = L_{k+2}(x) - L_{k}(x),\quad 0\le k\le N-2,
  \end{equation*}
and write
\begin{equation}\label{2.5}
\begin{split}
& u_{N} = \sum^{N-2}_{k,j = 0}\hat{u}_{kj}\phi_{k}(x)\phi_{j}(y),    \quad  U = (\hat{u}_{kj})_{k, j = 0, 1, \cdots, N-2};\\
& a_{kj} = \int_{I}\phi'_{j}(x)\phi'_{k}(x)dx,   \quad\quad  A = (a_{kj})_{k, j = 0, 1, \cdots, N-2};\\
& b_{kj} = \int_{I}\phi_{j}(x)\phi_{k}(x)dx,   \quad\quad  B = (b_{kj})_{k, j = 0, 1, \cdots, N-2};\\
& f_{kj} = (I_{N}f, \phi_{k}(x)\phi_{j}(y)),    \quad\;  g_{kj} = (I_{N}F(u_N), \phi_{k}(x)\phi_{j}(y))_{k, j = 0, 1, \cdots, N-2}.
\end{split}
\end{equation}
We obtain the nonlinear system
\begin{equation}\label{linsystemA}
(A \otimes B + B \otimes A^{\top})\bs{u} = \bs{f}-\bs {\hat{g}}(\bs u),
\end{equation}
where
\begin{equation*}
\begin{split}
& \bs{u} = (\hat{u}_{00}, \hat{u}_{10}, \cdots, \hat{u}_{q0}, \hat{u}_{01}, \cdots, \hat{u}_{q1}, \cdots, \hat{u}_{0q}, \cdots, \hat{u}_{qq})^{\top};\\
&\bs{f} = (f_{00}, f_{10}, \cdots, f_{q0}, f_{01}, \cdots, f_{q1}, \cdots, f_{0q}, \cdots, f_{qq})^{\top};\\
&\bs {\hat{g}}(\bs u) = (g_{00}, g_{10}, \cdots, g_{q0}, g_{01}, \cdots, g_{q1}, \cdots, g_{0q}, \cdots, g_{qq})^{\top}.
\end{split}
\end{equation*}
Here $\otimes$ denotes the Kronecker product, i.e. $A \otimes B = (Ab_{ij})_{i,j = 0, 1, \cdots, q}$ with $q = N-2$.

The nonlinear system here in \eqref{linsystemA} will be solved by the trust-region method.
Essentially, we only need to evaluate $\bs {\hat{g}}(\bs u)$ for given $\bs u$ that can be implemented efficiently by the pseduo-spectral technique described in   \cite[Ch.\! 4]{shen2011spectral}.
\begin{rem}\label{PseuSP} \emph{For the problem \eqref{New1.1} with high-order derivatives, it is advantageous to use spectral methods, which we shall describe
in Section \ref{sect3}.} \qed
\end{rem}

\subsection{Nonlinear least-squares and trust-region method}\label{TRM}
Formally, we write \eqref{linsystemA} as the following system of nonlinear equations:
\begin{equation}\label{4.1}
{\bs F}(\bs x) = \bs 0,  \quad  {\bs F} =\big({F}_{1}, \; {F}_{2}, \; \cdots, \; {F}_{n}\big)^{\top},
\end{equation}
with the unknown vector $  {\bs x} = (x_{1},\, \cdots,\, x_{n})^{\top}.$ We reformulate the zero-finding problem as the optimisation problem: 
\begin{equation}
\min_{\bs x \in {\mathbb R}^{n}}Q(\textbf{\textit{x}}),\quad Q(\bs x):= \frac{1}{2}\big\|\bs{F}(\bs x)\big\|^{2} = \frac{1}{2}\sum_{i = 1}^{n}{F}^2_{i}(\textbf{\textit{x}}), \label{4.2}
\end{equation}
where $\|\cdot\|$ is the vector 2-norm.
Denote  the gradient  and  Hessian matrices  by
\begin{equation*}
\begin{split}
& {\bs g}(\textbf{\textit{x}}) := \nabla Q(\textbf{\textit{x}}) = \bs{J}^{\top}(\bs{x}) {\bs F(\bs {x})} ,\\
& \bs G (\textbf{\textit{x}}) := \nabla^2 Q(\textbf{\textit{x}}) = {\bs J}^\top(\bs{x}) \bs J(\bs x) + \bs S(\bs x),
\end{split}
\end{equation*}
where
\begin{equation*}
\bs J(\bs x) = {\bs F}'(\textbf{\textit{x}}) = (\nabla {F}_{1}(\textbf{\textit{x}}), \nabla {F}_{2}(\textbf{\textit{x}}), \cdots, \nabla F_{n}(\textbf{\textit{x}}))^\top,  \quad {\bs S(\textbf{\textit{x}})} = \sum_{i = 1}^{n}F_{i}(\textbf{\textit{x}})\nabla^2 F_{i}(\textbf{\textit{x}}).
\end{equation*}

We resort to the trust-region method to solve the nonlinear least-squares minimization problem \eqref{4.2}. For this purpose, we introduce a region around the current best solution, called the trust region (say of radius $h_{k}$),  and approximate the objective function by a quadratic form which boils down to solving a sequence of trust-region subproblems:
\begin{equation}\label{3.2.1}
\begin{split}
& \min_{\bs s\in {\mathbb B}_{h_k}} q^{(k)}({\bs s}) := Q(\textbf{\textit{x}}^{(k)}) + {\bs g}(\textbf{\textit{x}}^{(k)})^{\top}{\bs s} + \frac{1}{2}{\bs s}^{\top}{\bs G}(\textbf{\textit{x}}^{(k)}){\bs s},\quad
\end{split}
\end{equation}
where the trust region  ${\mathbb B}_{h_k}:=\{\bs s\in {\mathbb R}^n\,:\,\|\bs s\|\le h_k\},$ and the matrices
${\bs g}(\textbf{\textit{x}}^{(k)})$ and ${\bs G}(\textbf{\textit{x}}^{(k)})$ are the gradient and  Hessian matrices at current point $\textbf{\textit{x}}^{(k)}$, respectively. Let ${\bs s}_{k}$ be the minimizer of $q^{(k)}({\bs s})$ in the trust-region of radius $h_{k}$. Then we update $\textbf{\textit{x}}^{(k+1)} = \textbf{\textit{x}}^{(k)} + {\bs s}_{k}$.  It is critical to choose a proper radius $h_k.$ In general, when there is good agreement between  $q^{(k)}({\bs s}_{k})$ and the objective function value $Q(\textbf{\textit{x}}^{(k+1)})$, one should choose $h_{k}$ as large as possible.  More precisely, assuming that
$q^{(k)}(0) \not= q^{(k)}({\bs s}_{k})$ (otherwise, ${\bs s}_{k}$ is a minimizer),  we define
\begin{equation}\label{rkval}
r_{k} = \frac{Q(\textbf{\textit{x}}^{(k)}) - Q(\textbf{\textit{x}}^{(k+1)})}{q^{(k)}(\bs 0) - q^{(k)}({\bs s}_{k})}.
\end{equation}
The ratio $r_{k}$ is an indicator  for the expansion and contraction the trust regions.   If $r_{k}$ is close to 1, it means there is good agreement, so we can expand the trust-region for the next iteration; if $r_{k}$ is close to zero or negative, we should shrink the trust-region radius; otherwise, we do not alter the trust-region.

%

We remark that as an efficient numerical optimization method for solving nonlinear programming (NLP) problems, the TRM enjoys the desirable global convergence with a local superlinear rate of convergence as follows.
\begin{thm}[see \cite{sun2006optimization}]\label{thm26}  Assume that
\begin{itemize}
 \item[(i)] the function $Q(\textbf{\textit{x}})$ is bounded below on the level set
\begin{equation}
S := \{\textbf{\textit{x}}\in R^{n}\; :\; Q(\textbf{\textit{x}}) \leq Q(\textbf{\textit{x}}^{(0)})\}, \quad \forall\, \textbf{\textit{x}}^{(0)}\in \mathbb R^n,
\end{equation}
and is Lipschitz continuously differentiable in $S;$
\item[(ii)] the Hessian matrixes $G(\textbf{\textit{x}}^{(k)})$ are uniformly bounded in 2-norm, i.e., $\|G(\textbf{\textit{x}}^{(k)})\| \leq \beta$ for any $ k$ and some  $\beta>0$.
\end{itemize}
If ${\bs g}(\textbf{\textit{x}}^{(k)}) \neq {\bs 0}$, then
\begin{equation}
\lim_{k \to \infty}\inf\|{\bs g}(\textbf{\textit{x}}^{(k)})\| = 0.
\end{equation}
Moreover,  if ${\bs g}(\textbf{\textit{x}}^{*}) = {\bs 0}$, and ${\bs G}(\textbf{\textit{x}}^{*})$ is positive definite,  then  the convergence rate of the TRM is quadratic.
\label{Nth1}
\end{thm}
\begin{rem}\label{remA} \emph{When $k$ is large enough, the TRM becomes the Newton iteration. As a result, it has the same convergence rate as the Newtonian method.} \qed
\end{rem}

\begin{rem} \emph{In practice,  the gradient and Hessian matrices might be appropriately approximated by some numerical means.
  We refer to Zhang et al. \cite{Hongchaozhang} for such derivative-free methods for \eqref{4.2} with
   $\bs F$ being   twice continuously differentiable, but none of their first-order or second-order derivatives being explicitly
available.} \qed
\end{rem}

\subsection{Deflation technique for multiple  solutions} 
\begin{figure}[H]
\begin{centering}
\includegraphics[width=7cm,height=4cm]{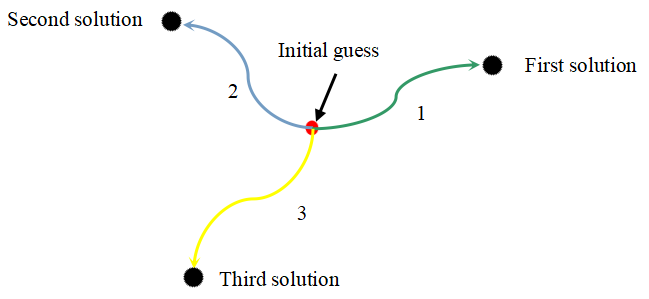}
\medskip
\caption{\small A schematic illustration of the deflation process integrated  with the LSTR method described in Subsection \ref{TRM}. Starting from the initial guess (marked by {\bf red dot}), we obtain the first solution (along {\bf Path 1}) via LSTR method. We then deflate this solution away by the deflation technique and update the objective function accordingly. As a result, we can find the second solution  (along {\bf Path 2}) via LSTR method with the same initial condition. Repeating this process, we search for other new solutions along different paths (e.g. {\bf Path 3)}.}\label{Fg3.1}
\end{centering}
\end{figure}

The critical issue in finding multiple solutions is to
 get away from one identified  solution and search for  a new solution based on ideally  the same initial input.  Motivated  by Farrell  et al \cite{farrell2015deflation}, we introduce an effective technique rooted in the notion of deflation. As illustrated  in Figure  \ref{Fg3.1}, suppose that the first (approximate) solution $\bs x_1^*$ of ${\bs F}(\bs x)=\bs 0$  in \eqref{4.1} via the LSTR method, and
  we intend to find a second (approximate) root $\bs x_2^*$ of  ${\bs F}.$ For this purpose, we consider the modified zero-finding problem $\widehat{\bs F}(\bs x)=\bs 0$ where $\widehat{\bs F}$ is constructed in a way to deflate the first solution $\bs x_1^*,$ and allow for  finding $\bs x_2^*$ by the LSTR method even with the same initial data.  We then repeat this deflation process to find more solutions.


We next introduce  the  deflation operator, and highlight the essence of the deflation technique
\cite{0Nonlinear,1971Deflation,farrell2015deflation}.

\vskip 4pt
\noindent\textbf{Definition of deflation operator.}
\emph{Let ${\mathbb X}, {\mathbb Y},$ and ${\mathbb Z}$ be Banach spaces, and $\widehat{\mathbb X}$ be an open subset of ${\mathbb X}$. Let ${\bs F}: {\mathbb X} \to {\mathbb Y}$ be a Fr\'{e}chet differentiable operator with the Fr\'{e}chet derivative ${\bs F'}$.
Suppose that  ${\bs r} \in \widehat{\mathbb X}$ is a root of ${\bs F},$ i.e.,  ${\bs F}(\bs r)=\bs 0,$ such that  ${\bs F}'({\bs r})$ is nonsingular.  For each ${\bs x}\in \widehat{\mathbb X} \setminus \{{\bs r}\},$ let
$\Psi({\bs x}; {\bs r}): {\mathbb Y} \to  {\mathbb Z}$ be an invertible linear operator,
and define
\begin{equation}
\widehat{{\bs F}}({\bs x}) = \Psi({\bs x}; {\bs r}){\bs F}({\bs x}).\label{New4.2.2}
\end{equation}
  We call  $\Psi$  a deflation operator, if  for any sequence $\{{\bs x}_{i}\}\subseteq \widehat{\mathbb X}\setminus \{{\bs r}\}$ converging to ${\bs r}$,  we have}
\begin{equation}
\lim_{i \to \infty}\inf \|\widehat{{\bs F}}({\bs x}_{i})\|_{\mathbb Z} > 0.  \label{New4.2.1}
\end{equation}
It is seen from (\ref{New4.2.2})-(\ref{New4.2.1}) that the role of the deflation operator is to exclude and deflate the root $\bs r,$ so
$\widehat{{\bs F}}({\bs x})$ has common roots as ${{\bs F}}({\bs x})$ except for $\bs r.$

We next present some sufficient conditions for  constructing suitable deflation operators.
\begin{thm}[see \cite{farrell2015deflation}]\label{thm22} Let ${\bs F}:$ $\widehat{\mathbb X} \subset {\mathbb  X} \to  {\mathbb  Y}$ be a Fr\'{e}chet differentiable operator. Suppose that the linear operator $\Psi({\bs x}; {\bs r}): {\mathbb Y} \to {\mathbb Z}$ has the property that for each ${\bs r} \in \widehat{\mathbb X}$, and any sequence ${\bs x}_{i} \rightarrow {\bs r}$, ${\bs x}_{i}\in \widehat{\mathbb X} \setminus \{{\bs r}\}$,
\begin{equation}
\text{if}\quad  \|{\bs x}_{i}-{\bs r}\|\Psi({\bs x}_{i}; {\bs r}) {\bs F}({\bs x}_{i}) \to 0  \quad \text{implies} \quad  \; {\bs F}({\bs x}_{i}) \to 0,\quad i\to \infty, \label{4.2.2}
\end{equation}
 then  $\Psi$ defines a deflation operator.\label{th1}
\end{thm}\vspace{-6pt}
Based on this rule,  several deflation operators have been suggested in \cite{farrell2015deflation}. In what follows, we adopt the shifted deflation operator:
\begin{equation}
\Psi({\bs x}; {\bs r}) = \frac{\mathbb I}{\|{\bs x}-{\bs r}\|^{2}} + {\mathbb I},
\end{equation}
where ${\mathbb I}$ is the identity operator on ${\mathbb Y}$. Suppose that $\bs r_1$ is a solution of the original nonlinear system $\bs F(\bs x)=\bs 0.$
Then we solve the deflated nonlinear system:  ${\bs F}_1(\bs x):=\Psi({\bs x}; {\bs r_1}) {\bs F}(\bs x)=\bs 0$ to obtain the second solution $\bs r_2.$ Repeating this process, we search for the $(k+1)$-th solution by solving the nonlinear system:
\begin{equation}\label{NNNew2.3.1}
{\bs F}_k(\bs x): = {\bs F}(\bs x) \prod \limits_{i=1}^{k}\Psi(\bs x; {\bs r_i}) =\bs 0.
\end{equation}

\subsection{Summary of the spectral LSTR-Deflation algorithm}\label{sub24}
In summary, we have the following algorithm for finding multiple solutions.

\smallskip
\begin{center}
\begin{tabular}{l}
\toprule[0.5pt]
Spectral LSTR-Deflation Algorithm for Computing Multiple Solutions\label{3.15}\\\toprule[0.5pt]
\quad \textbf{Input:} Given $\textbf{\textit{x}}^{(0)}$, $\epsilon > 0,  0<\delta_1 < \delta_2 < 1, 0 < \tau_1 < 1 < \tau_2,$ $\bs F (\bs x)=\bs 0,$
and $h_{0} = \|\bs g_{0}\|$\\ 
\quad \textbf{Input:} Initial solution set  $S \leftarrow$ the empty set $\Phi$\\
\quad \textbf{Output:}  $S$ \\
\,\,\, 1 {\color{red}{while}} the desired multiple solutions are not found {\color{red}{do}}\\
\,\,\, 2 \quad {\color{blue}{while}} first-order optimality threshold or failure criterion does not meet {\color{blue}{do}}\\
\,\,\, 3  \quad\quad Compute ${\bs g}_{k}$ and ${\bs G}_{k}$;\\
\,\,\, 4  \quad\quad If $\|{\bs g}_{k}\| \leq \epsilon$ and $|Q(\textbf{\textit{x}}^{(k)})| \leq \epsilon$, stop;\\
\,\,\, 5  \quad\quad Approximately solve the subproblem (\ref{3.2.1}) for ${\bs s}_{k}$;\\
\,\,\, 6  \quad\quad Compute $Q(\textbf{\textit{x}}^{(k)} + {\bs s}_{k})$ and $r_{k}$. If $r_{k}\geq \delta _1,$ then $\textbf{\textit{x}}^{(k+1)} = \textbf{\textit{x}}^{(k)} + {\bs s}_{k}$; Otherwise,\\
\,\,\,   \quad\quad\quad set $\textbf{\textit{x}}^{(k+1)} = \textbf{\textit{x}}^{(k)}$;\\
\,\,\, 7  \quad\quad If $r_{k} < \delta_1$, then $h_{k+1} = \tau_1 h_{k}$;\\
\,\,\,   \quad\quad\quad If $r_{k} > \delta_2$ and $\|{\bs s}_{k}\|$ = $h_{k}$, then $h_{k+1} = \tau_2 h_{k}$;\\
\,\,\,   \quad\quad\quad Otherwise, set $h_{k+1} = h_{k}$;\\
\,\,\, 8  \quad\quad Generate ${\bs G}_k$, update $q^{(k)}$, set $k := k + 1$, go to step 3;\\
\,\,\, 9   \quad {\color{blue}{end}}\\
\, 10  \quad \textbf{if} Convergence threshold met \textbf{then}\\
\, 11  \quad\quad\quad Add converged solution $\tilde{\textbf{\textit{x}}}$ to $S$\\
\, 12   \quad\quad\quad Set $\textbf{\textit{F}}$(\textbf{\textit{x}}) $\leftarrow$ $\Psi(\textbf{\textit{x}}; \tilde{\textbf{\textit{x}}})$\textbf{\textit{F}}(\textbf{\textit{x}}) (i.e. \textbf{Deflation}), and go back to step 2.\\
\, 13   \quad \textbf{end}\\
\, 14    {\color{red}{end}}\\
\, 15  return $S$\\
\toprule[0.5pt]
\end{tabular}
\end{center}
In the above,  we have ${\bs g}_{k} = {\bs g}(\textit{\textbf{x}}^{(k)})$ and ${\bs G}_{k} = {\bs G}(\textit{\textbf{x}}^{(k)})$. For the trust-region subproblem (see Line 5), an efficient implementation for its solution is  the so-called dogleg method (see \cite{sun2006optimization, 0a}) with the process:  
\begin{equation*}
q^{(k)}({\bs x}^{(k)} - l_{k} {\bs g}_{k}) = Q({\bs x}^{(k)}) - l_{k}\|{\bs g}_{k}\|^{2} + \frac{1}{2}l^2_{k} {\bs g}^{\top}_{k}{\bs G}_{k}{\bs g}_{k},
\end{equation*}
where based on the exact line search, we have
\begin{equation*}
l_{k} = \frac{\|{\bs g}_{k}\|^{2}}{{\bs g}^{\top}_{k}{\bs G}_{k}{\bs g}_{k}}.
\end{equation*}
The corresponding step along the steepest descent direction is
\begin{equation*}
{\hat{\bs s}}_{k} = - l_{k}{\bs g}_{k} = -\frac{{\bs g}^{\top}_{k}{\bs g}_{k}}{{\bs g}^{\top}_{k}{\bs G}_{k}{\bs g}_{k}}{\bs g}_{k},
\end{equation*}
and  the Newton step is
\begin{equation*}
\tilde{{\bs s}}_{k} = -{\bs G}^{-1}_{k}{\bs g}_{k}.
\end{equation*}
If $\|\hat{{\bs s}}_{k}\|_{2} = \|l_{k}{\bs g}_{k}\|_{2} \geq h_{k}$, we take
\begin{equation}
{\bs s}_{k} = -\frac{h_{k}}{\|{\bs g}_{k}\|}{\bs g}_{k} \quad \textrm{and} \quad  {\bs x}^{(k+1)} = {\bs x}^{(k)} + {\bs s}_{k}.\label{New43}
\end{equation}
If $\|\hat{{\bs s}}_{k}\| < h_{k}$ and $\|\tilde{{\bs s}}_{k}\| > h_{k}$, we take
\begin{equation}
{\bs s}_{k}(\lambda) = \hat{{\bs s}}_{k} + \lambda(\tilde{{\bs s}}_{k} - \hat{{\bs s}}_{k}),\quad  0 \leq \lambda \leq 1.\label{New42}
\end{equation}
As a result, we have
\begin{equation*}
{\bs x}^{(k+1)} = {\bs x}^{(k)} + {\bs s}_{k}(\lambda) = {\bs x}^{(k)} + \hat{{\bs s}}_{k} + \lambda(\tilde{{\bs s}}_{k} - \hat{{\bs s}}_{k}),
\end{equation*}
where $\lambda$ is obtained by solving the equation
\begin{equation*}
\|\hat{{\bs s}}_{k} + \lambda(\tilde{{\bs s}}_{k} - \hat{{\bs s}}_{k})\| = h_{k}.
\end{equation*}
Otherwise, we set
\begin{equation}
{\bs s}_{k} = \tilde{{\bs s}}_{k} = -{\bs G}^{-1}_{k}{\bs g}_{k}. \label{New41}
\end{equation}
In summary, combining (\ref{New43}), (\ref{New42}) and (\ref{New41}) yields
\begin{equation}
{\bs s}_{k}=\begin{dcases}
  -\frac{{\bs g}_{k}}{\|{\bs g}_{k}\|} h_{k}, &  \textrm{if}\;\; \|\hat{{\bs s}}_{k}\|_{2} \geq h_{k}, \\
\hat{{\bs s}}_{k} + \lambda(\tilde{{\bs s}}_{k} - \hat{{\bs s}}_{k}), & \textrm{if}\;\; \|\hat{{\bs s}}_{k} \|_{2} < h_{k}\;\;
 \textrm{and}\;\;  \|\tilde{{\bs s}}_{k}\|_{2} > h_{k},  \\
 -{\bs G}^{-1}_{k}{\bs g}_{k}, &  \textrm{if}\;\; \|\hat{{\bs s}}_{k}\|_{2} < h_{k}\;\;  \textrm{and}\;\;  \|\tilde{{\bs s}}_{k}\|_{2} \leq h_{k}.
\end{dcases}
\end{equation}

Some  additional remarks about the implementation of the Spectral LSTR-Deflation Algorithm are listed as follows.
\begin{itemize}
  \item  As highlighted previously, the choice of the trust-region radius $h_k$ in the LSTR method is important.
   The ratio $r_{k}$ in \eqref{rkval} is used as the criterion for expansion or contraction of the region.
  \smallskip
  \item For the original nonlinear algebraic system ${\bs F}$ (see Line 12 of the Algorithm), the solution $\tilde{{\bs x}}$ is found from initial guess ${\bs x}^{(0)}$ with the Fr\'{e}chet derivative ${\bs F}'(\tilde{{\bs x}})$ being nonsingular. While here the deflation operator $\Psi$ plays a key role in finding multiple solutions other than $\tilde{{\bs x}}$. Let $\hat{{\bs F}}({\bs x}) = \Psi({\bs x}; \tilde{{\bs x}}){\bs F}({\bs x})$, then this deflated problem $\widehat{{\bs F}}({\bs x})$ satisfies two properties: (i) The preservation of solutions of $\hat{{\bs F}}$ should be hold, i.e. for ${\bs x} \neq \tilde{{\bs x}}$, $\widehat{{\bs F}}({\bs x}) = 0$ iff ${\bs F}({\bs x}) = 0$; (ii) In view of  Theorem \ref{th1}, the iterative algorithm applied to $\widehat{{\bs F}}({\bs x})$ will not find $\tilde{{\bs x}}$ again, when
\begin{equation*}
\lim_{{\bs x} \to \tilde{{\bs x}}}\inf \|\widehat{{\bs F}}({\bs x})\| > 0,
\end{equation*}
i.e., the deflated function $\widehat{{\bs F}}({\bs x})$ does not vanish at the known solution $\tilde{{\bs x}}$.
\smallskip
  \item In general, based on the deflation technique, the same initial guess allows us to find new solutions in most cases. However, it is noteworthy that the deflation procedure may  diverge with the same initial guess in some cases. We find out that choosing a random perturbation of an obtained solution as the initial guess for the deflated problems turns out sufficient for finding all other solutions (see Section \ref{sect3}).
\smallskip
\item Although the original and deflated nonlinear systems can be solved  by other methods, e.g., Newton iteration, the trust-region method is deemed as the method of choice for its  flexibility in initial inputs and   efficiency in dealing with the increasingly complicated systems in the course of repeating the deflation technique. We shall testify this in the forthcoming section.
\end{itemize}

\smallskip
\section{Applications to differential equations with multiple solutions}\label{sect3}
\setcounter{equation}{0}
\setcounter{lmm}{0}
\setcounter{thm}{0}
In this section, we apply the spectral LSTR-Deflation algorithm to various differential equations with multiple solutions and compare its performance in efficiency and accuracy with some related methods.
In the numerical experiments, we choose  $\delta_1 = 0.25, \delta_2 = 0.75, \tau_1 = 0.5$, $\tau_2 = 2$ and $\epsilon = 10^{-13}$ in  the algorithm described in Subsection \ref{sub24}. The codes  are carried out on a server: Intel(R) Core(TM) i7-7500U (2.90 GHz) and 120GB RAM
via  MATLAB (version R2015b).

\subsection{ODE Examples}
\subsubsection{Viscous flow in a porous channel with expanding or contracting walls} When concerning systemic circulation in blood circulation, the blood in the left ventricle is being forced into the aorta by systole and the mitral valve between left ventricle and left atrium is closed. At this juncture the left ventricle forms a vessel with one end closed. Meanwhile, the mass transfer of the vessel between inside and outside can be achieved by the seepage across permeable wall of the vessel. The mathematical model  established by Majdalani
\cite{dauenhauer2003exact} is as follows 
\begin{equation}
u^{(4)}(y) + \alpha (yu''' (y)+ 3u''(y)) + Re\big(u(y)u'''(y) - u'(y)u''(y)\big) = 0, \;\;\; y\in (0,1), \label{1.1}
\end{equation}
supplemented with the  boundary conditions
\begin{equation}
u(0) = 0, \quad  u''(0) = 0, \quad u(1) = 1,  \quad  u'(1) = 0,  \label{1.2}
\end{equation}
where $Re$ is called as the cross-flow Reynolds number ($Re > 0$ for injection and $Re < 0$ for suction), $\alpha$ is the wall expansion ratio and the prime denotes differentiation with respect to $x$. Here $u(y)$ and $u'(y)$ represent the normal and streamwise velocities in a porous channel,
respectively. this problem admits multiple solutions for $Re < -12.165$ and any $\alpha$, but so far there are only three solutions identified in literature \cite{robinson1976existence, xu2010homotopy}.


\medskip

\noindent{\bf (i)\, Legendre spectral Petrov-Galerkin discretisation.}\;
As shown in  \cite{shen2011spectral}, the spectral method is advantageous for high-order problems.  In view of the non-symmetric boundary conditions \eqref{1.2}, we employ the Petrov-Galerkin method using compact combinations of Legendre polynomials that meet the boundary conditions and the dual counterparts.
For convenience, we make a change of variable: $y=(1+x)/2,$ and
convert (\ref{1.1})-(\ref{1.2}) into
\begin{equation} \label{2.2}
\begin{dcases}
 \breve u^{(4)} + \delta(1+x){\breve u}''' + \beta {\breve u}'' + \gamma F({\breve u}) = 0, \quad x\in (-1,1),\\
 {\breve u}(-1) = {\breve u}''(-1) = 0,   \quad {\breve u}(1) = 1,    \quad   {\breve u}'(1) = 0,
\end{dcases}
\end{equation}
where ${\breve u}(x) = u(y)$, $\delta = \frac{\alpha}{4}$, $\beta = \frac{3}{2}\alpha$, $\gamma = \frac{Re}{2}$, and the nonlinear term:
\begin{equation}\label{34form}
F({\breve u}) := {\breve u}{\breve u}''' - {\breve u}'{\breve u}'' = \frac{1}{2}({\breve u}^2)''' - 2(({\breve u}')^2)'.
\end{equation}
To deal with the non-homogeneous boundary conditions, we introduce
\begin{equation*}
p(x) = \frac{3}{4}(1+x) - \frac{1}{16}(1+x)^3,\;\;\; {\rm s.t.,}\;\;\; p(-1) = p''(-1) = p'(1) = 0,  \;\;  p(1) = 1.
\end{equation*}
Set ${\breve u}=v+p,$ and note
\begin{equation*}
F(v+p) = F(v) + pv''' - p'v'' - p''v' + p'''v + pp''' - p'p''.
\end{equation*}
Then direct calculation leads to
\begin{equation}\label{u4eqn00}
\begin{cases}
 v^{(4)}+A(x) v^{\prime \prime \prime}+B(x) v^{\prime \prime}+C(x) v'-\dfrac 3 8 \gamma v+\gamma F(v) = g(x),\quad x\in (-1,1),\\[2pt]
 v(\pm 1)=v'(1)=v''(-1)=0,
\end{cases}
\end{equation}
where
\begin{equation}\label{u4eqnf1}
\begin{split}
& A(x)=\delta(1+x)+\gamma p(x)=\Big(\frac {3\gamma} 4 +\delta\Big)(1+x)- \frac{\gamma}{16}(1+x)^3,\\[4pt]
& B(x)=\beta-\gamma p'(x)=\beta-\frac {3\gamma} 4  +\frac {3\gamma}{16}(1+x)^2,\quad C(x)=\frac{3\gamma} 8(1+x),  \\[4pt]
&g(x)=\frac {3} 8 (\delta+\beta)(1+x) -\gamma pp'''+\gamma p'p''=  \frac {3} 8 (\delta+\beta)(1+x)+\frac {3\gamma}{64} (1+x)^3.
\end{split}
\end{equation}

Let $\mathbb P_N$ be the set of polynomials of degree $\le N$, and define
\begin{equation}\label{VnWn}
\begin{split}
&V_N:=\big\{\phi\in \mathbb P_N:  \phi(\pm 1)=\phi'(1)=\phi''(-1)=0
\big\},\\
&W_N:=\big\{\psi\in \mathbb P_N:  \psi(\pm 1)=\psi'(\pm 1)=0
\big\}.
\end{split}
\end{equation}
The  Legendre spectral Petrov-Galerkin scheme for  \eqref{u4eqn00}-\eqref{u4eqnf1} is to find $v_N\in V_N$  such that
\begin{equation}\label{scheme}
\begin{split}
{\mathcal B}(v_N, \psi)+ {\mathcal N}(v_N, \psi)=(g,\psi),\quad \forall \psi\in W_N,
\end{split}
\end{equation}
where the bilinear form and the approximation to the nonlinear term are
\begin{equation}\label{scheme2}
\begin{split}
& {\mathcal B}(v_N, \psi):=(v_N'', \psi'')+(v_N',(A\psi)'')- (v_N',(B\psi)') +(Cv_N', \psi)-\frac {3\gamma} 8 (v_N,\psi), \\
&{\mathcal N}(v_N, \psi):= \frac{\gamma}2 \big( [I_N v_N^2]',\psi''\big) +2\gamma \big(I_N [v_N']^2, \psi'),
\end{split}
\end{equation}
with $I_N$ being the Legendre-Gauss-Lobatto interpolation operator. Here, we use \eqref{34form} to formulate ${\mathcal N}(v_N, \psi).$

%

 To construct basis functions, we now recall some  properties of Legendre polynomials (cf.\!  \cite{shen2011spectral}): 
\begin{equation}\label{legendA}
\begin{split}
& (2n+1)L_n(x)=L_{n+1}'(x)-L_{n-1}'(x)=J_{n-1}'(x),\\
& (1-x^2)L_n'(x)=\frac{n(n+1)}{2n+1}\big(L_{n-1}(x)-L_{n+1}(x)\big)=-\frac{n(n+1)}{2n+1} J_{n-1}(x),
\end{split}
\end{equation}
and
\begin{equation}\label{legendB}
L_n(\pm 1)=(\pm 1)^n,\quad  L'_n(\pm 1)=\frac12(\pm 1)^{n-1}n(n+1),
\end{equation}
where  we denoted
\begin{equation}
J_n(x)=L_{n+2}(x)-L_{n}(x),\quad n\ge 0.
\end{equation}
Direct calculation yields
\begin{equation}\label{Jvalues}
\begin{split}
& J_n'(1)= (2n+3) L_{n+1}(1)=2n+3, \\[4pt]
& J_n''(-1)= (2n+3) L_{n+1}'(-1)= \frac12(-1)^{n}(n+1)(n+2) (2n+3).
\end{split}
\end{equation}
Now, we construct the basis functions for $V_N$ using the  compact combination: 
\begin{equation}\label{legendB0}
\varphi_k(x)=J_k(x)+ a_k J_{k+1}(x)+ b_k J_{k+2}(x), \quad 0\le k\le N-4,
\end{equation}
where $\{a_k, b_k\}$ are determined by the homogenous boundary conditions in $V_N.$  In view of the second formula in  \eqref{legendA}, we have  $\varphi_k(\pm 1)=0.$ From \eqref{legendA}-\eqref{legendB}   and  $\varphi_k'(1)=\varphi_k''(-1)=0,$  we can obtain
\begin{equation*}
J_{k+1}'(1)\, a_k+ J_{k+2}'(1)\, b_k=- J_{k}'(1),\quad
J_{k+1}''(-1)\, a_k+ J_{k+2}''(-1)\, b_k=- J_{k}''(-1),
\end{equation*}
which, together with \eqref{Jvalues},  implies
\begin{equation}\label{legendC}
\begin{split}
& a_k=\frac{J_{k+2}'(1)J_{k}''(-1)- J_{k}'(1)J_{k+2}''(-1)}{J_{k+1}'(1)J_{k+2}''(-1)- J_{k+2}'(1)J_{k+1}''(-1) }=-\frac{2k+3} {(k+3)^2},\\[6pt]
&  b_k=\frac{J_{k}'(1)J_{k+1}''(-1)- J_{k+1}'(1)J_{k}''(-1)}{J_{k+1}'(1)J_{k+2}''(-1)- J_{k+2}'(1)J_{k+1}''(-1) }=-\frac{(k+2)^2(2k+3)}{(k+3)^2(2k+7)}.
\end{split}
\end{equation}
Similarly, we can also construct the basis functions for $W_N:$
\begin{equation}\label{legendC0}
\psi_k(x)=J_k(x)+  c_k J_{k+2}(x), \quad c_k=-\frac {2k+3}{2k+7}.
\end{equation}

We expand the numerical solution as
\begin{equation}\label{unumerical }
u_N(x):=u_N(x;\alpha,Re)=p(x)+v_N(x)=p(x)+\sum_{k=0}^{N-4} \tilde v_k \varphi_k(x).
\end{equation}
Under these basis functions, we can easily show that the matrix of the linear part  ${\mathcal B}(\varphi_k, \psi_j)$ in
\eqref{scheme} is sparse with a finite bandwidth and its entries can be evaluated explicitly.  It is evident that the righthand sided vector
$\{(g,\psi_j)\}_{j=0}^{N-4}$ can be computed exactly as well.
%
Then we solve the nonlinear system by the LSTR-Deflation technique summarized in Subsection \ref{sub24},
where we choose the initial guess
\begin{equation}\label{N1.1}
\tilde v_k^{(0)}=1,\quad 0\le k\le N-4,
\end{equation}
denoted by $\textrm{IG} = \textrm{ones}(N-3).$ Note that in the iterations, the nonlinear part $\{{\mathcal N}(v_N, \psi_j)\}_{j=0}^{N-4}$
can be treated efficiently by using the pseudospectral technique as in \cite[Ch.\! 4]{shen2011spectral}.

\medskip

\begin{figure}[!h]
\centering
\subfigure[$\alpha = 0, Re = -20$]{
\includegraphics[width=6.5cm]{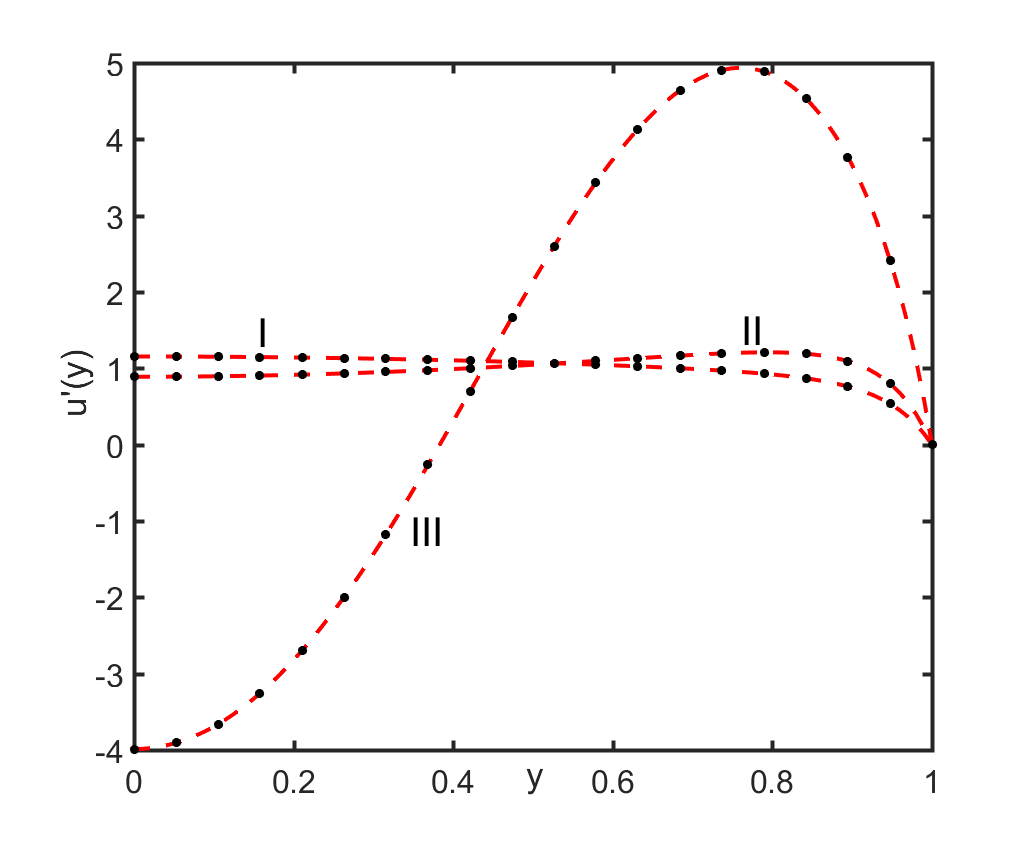}
    }
\quad
\subfigure[$\alpha = 8, Re = -40$]{
\includegraphics[width=6.5cm]{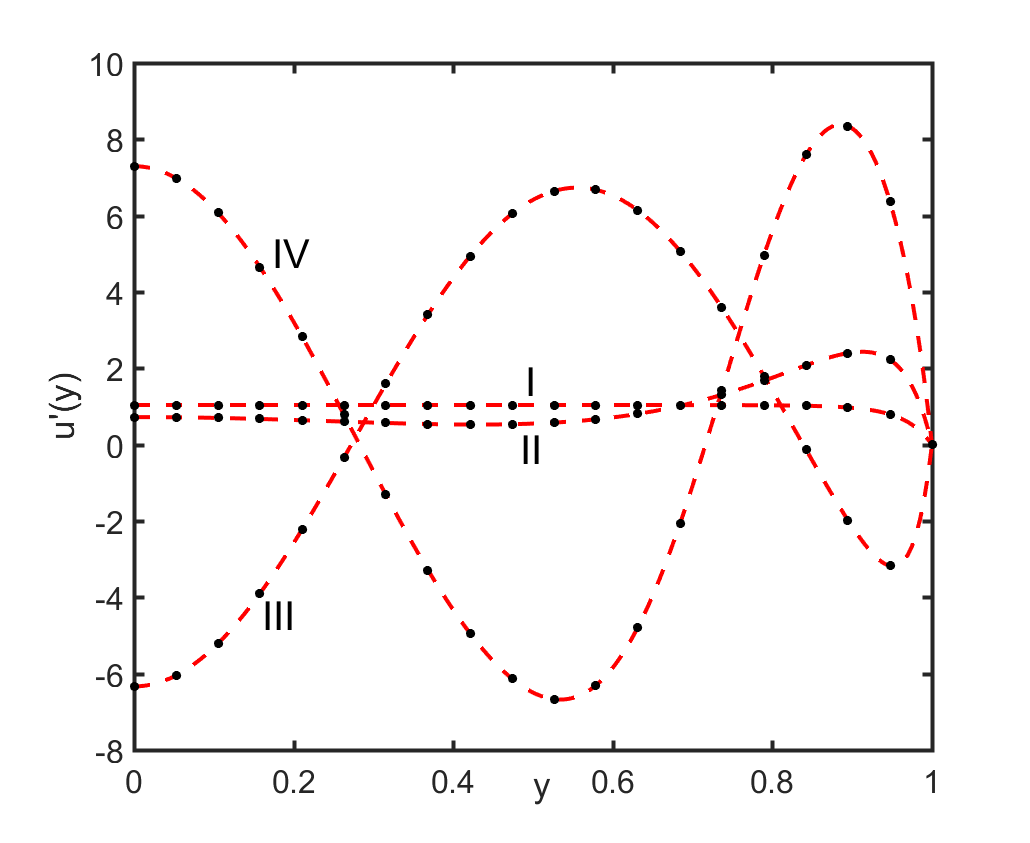}
}
\caption{Multiple solutions of (\ref{1.1})-(\ref{1.2}) with different ($\alpha, Re$).}\label{T3}
\end{figure}

\noindent{\bf (ii)\, Numerical results.}\;
In Table \ref{Table42}, we tabulate the outcomes of our algorithm for four pairs of parameters $(\alpha, Re)$ and the above initial guess with $N=18,$
where we list the number of deflations,  solutions (labeled by  {\rm I,II,$\cdots$}), the computational time and  iteration number of the LSTR method for each solution.
Some observations and highlights are as follows.
\begin{itemize}
\item[i)] It is noteworthy that in \cite{robinson1976existence, xu2010homotopy}, the plots of the solutions for
$(\alpha, Re)=(0, -20), (2,-40)$ using different methods were reported.  Here, we obtain the solutions in Figure \ref{T3} efficiently in seconds from the same initial guess
$\textrm{IG}$.
\item[ii)] We find a new solution  for $(\alpha, Re)$=(8, -40), see solution-IV in Figure \ref{T3}.
\item[iii)] The algorithm is quite accurate with a fast convergence  as shown in Table \ref{NTable41}, where we compare the errors in maximum norm with the numerical solution obtained with a relatively large $N.$
\end{itemize}

\begin{table}[!h]
\centering\small
\caption{\small  Performance of our spectral LSTR-Deflation method.}
\label{Table42} \small 
  \begin{tabular}{|c|cccc|c|cccc|}
	\hline 
	$(\alpha, Re)$ &Deflations\!\!  &Solutions\!\! & Time(s)\!\! &$n_{it}$   &$(\alpha, Re)$ & Deflations\!\!   &Solutions\!\! & Time(s)\!\! &$n_{it}$\\ 
	\hline
	\multirow{3}{*}{(0, -20)}  &\multirow{3}{*}{2} &$\textrm{I}$ & 3.4567 & 24 &\multirow{3}{*}{(2, -40)}   &\multirow{3}{*}{2} & $\textrm{I}$ & 3.2326 & 24 \\
	 ~ &~  & $\textrm{II}$ & 2.7378 & 17  &~ &~ &$\textrm{II}$ &5.2305 & 43 \\
	 ~ &~  &$\textrm{III}$ & 2.2060 & 12 &~ &~ & $\textrm{III}$  &4.4526 &34  \\
\hline
     \multirow{4}{*}{(-2, -40)}  &\multirow{4}{*}{2} & $\textrm{I}$ & 4.0644 & 31 & \multirow{4}{*}{(8, -40)}    &\multirow{4}{*}{3}  & $\textrm{I}$ & 3.0578 & 22\\
	 ~ &~ &$\textrm{II}$ & 4.1458 & 31 &~ &~ &$\textrm{II}$ & 4.1466 &35 \\
	 ~ &~ & $\textrm{III}$ & 4.4088 & 33 &~ &~ & $\textrm{III}$ & 4.1352 &34 \\
	 ~ &~ & ~ & ~ & ~  & ~ &~ & $\textrm{IV}$ & 5.1466 & 42\\
\toprule[0.5pt]
\end{tabular}
\end{table}

\begin{table}[!h]
\centering
\caption{\small \normalsize Accuracy of our spectral LSTR-Deflation method.}
\label{NTable41}\small
\begin{tabular}{|c c c c c|c c c c c c|}
	\hline
$(\alpha, Re)$ &$N$ &I  &II   &III &$(\alpha, Re)$ &$N$ &I  &II   &III  &IV \\\hline
\multirow{3}{*}{(0, -20)}  & 8   &7.71e-2   &1.06e-1   &2.36e-3 &\multirow{3}{*}{(2, -40)} & 8   & 3.15e-2   & 1.23e-2   & 4.25e-2  &-\\
~   & 16   & 6.32e-5   & 1.02e-5   & 3.16e-7 &~  & 16   & 3.67e-5   &7.21e-5   &3.81e-5 &-     \\
~   & 24   & 5.21e-9   & 1.78e-8   & 2.08e-10  &~  & 24   &4.31e-8   &9.23e-7   &3.71e-9  &- \\ \hline
\multirow{3}{*}{(-2, -40)}  &8  & 2.81e-2   &3.61e-1   &4.27e-1 &\multirow{3}{*}{(8, -40)}  & 8   &1.49e-1   & 3.07e-1   &2.89e-1  &9.31e-1\\
~   & 16   & 5.31e-5   & 2.39e-5   &3.86e-7 &~ & 16   &1.83e-5   &4.31e-5   &1.09e-5 &4.81e-6\\
~   & 24   & 3.89e-8   &1.84e-8   &4.61e-9 &~  & 24   & 3.91e-8   & 3.92e-9   &9.18e-9   &1.78e-8\\ \hline
\end{tabular}
\end{table}


As a comparison, we replace the LSTR algorithm for solving the nonlinear systems by the Newton iteration (as in \cite{farrell2015deflation} using Newton-Deflation technique). We compute  $\|{\bs F}(\bs x)\|$ in \eqref{4.2} for three choices of initial guesses  ${\rm IG}, 0.1*{\rm IG}$ and $\cos({\rm IG})$.
As shown in Table \ref{Table41}, the Newton iterations fail to converge for these initial inputs, while  $\|{\bs F}(\bs x)\|$ of the
LSTR method descends very fast. However, it is noteworthy that if we choose the initial inputs
to be a small perturbation of the numerical solution from our algorithm, the Newton-Deflation method indeed converges.

%
%
%

\begin{table}[!h]
\centering\small
\caption{\small A comparison between Newton iteration in \cite{farrell2015deflation} and the LSTR method.}
\label{Table41}
\begin{tabular}{|ccccc|cccc|}
	\hline
~ &\multicolumn{4}{c}{Newton iteration in \cite{farrell2015deflation}} &\multicolumn{4}{c|}{LSTR method}\\\cline{2-9}
$(\alpha, Re)$ &$n_{it}$ & ${\rm IG}$  & $0.1 *{\rm IG}$   & $\cos({\rm IG})$  &$n_{it}$  & ${\rm IG}$  & $0.1 *{\rm IG}$   & $\cos({\rm IG})$ \\  \hline
\multirow{5}{*}{(8, -40)}  & 1   & 4.97e5   & 2.02e3   & 3.13e6  &10   & 2.57e-3   & 3.72e-13   & 8.87e1\\
~   & 2   & 3.10e9   & 8.59e3   & 9.91e6   &15   & 6.50e-5   &-  & 2.45e-4     \\
~   & 3   & 3.39e10   & 5.07e4   & 2.17e9   &20   & 4.41e-13   &-   &4.98e-6   \\
~   & 4   & 2.64e13   & 2.16e5   & 3.03e11   &25   &-   &-   &6.91e-9  \\
~   & 5   & 1.99e15   & 1.47e5   & 5.66e13   &30   &-  &-   &1.86e-13 \\\hline
\end{tabular}
\end{table}

\subsubsection{Models with polynomial and exponential nonlinearities}
We consider the  model problem:
\begin{equation}\label{3.5}
u'' + \lambda F(u) = 0,  \quad x \in (0,1), \quad \lambda>0,
\end{equation}
with boundary conditions to be specified later.  Here, we focus on $F(u) = e^{u}, 1+u^{p}$ (where $p$ is a positive integer).

\medskip
\noindent\underline{\bf Case 1:}\;  The equation \eqref{3.5} with $F(u) = e^{u}$ is known as the  Bratu-Gelfand model  \cite{davis1960introduction}, which
together with the boundary conditions $u(0) = u(1) = 0,$ 
 has two solutions when  $0 < \lambda < \lambda^{*} (\lambda^{*} = 3.51383)$, but no solution for $\lambda \ge \lambda^{*}$. Here we assume that $0 < \lambda < \lambda^{*}$. The spectral-Galerkin discretisation follows   \cite[Ch.\! 4]{shen2011spectral}.


\begin{table}[!h]
\centering\small
\caption{\small Performance of the spectral LSTR-Deflation methods.}
\label{NewTable61}
\begin{tabular}{|cccccc|cccc|}
	\hline
~ &~ &\multicolumn{4}{c}{First solution} &\multicolumn{4}{c|}{Second solution}\\\cline{3-10}
\specialrule{0em}{0.1pt}{0.1pt}
$\lambda$ &$N$  &$n_{it}$ & Time(s)   &$L^{\infty}$-error &$\|{\bs F}({\bs x})\|_{\infty}$  &$n_{it}$ & Time(s) &$L^{\infty}$-error   &$\|{\bs F}({\bs x})\|_{\infty}$ \\\hline
\multirow{3}{*}{1} &8  &3  &0.0367 &2.6788e-7  &4.9021e-12    &20  &0.1623 &2.0258e-4  &8.1304e-12 \\
~ &16  &4  &0.0492 &2.2355e-9  &2.9015e-13                  &26  &0.1878 &1.3531e-4  &5.9026e-13 \\
~ &24  &6  &0.0512 &5.9229e-11  &1.8935e-12                  &28  &0.1886 &1.1749e-7   &3.9057e-14\\\hline
\multirow{3}{*}{2} &8  &4  &0.0501 &1.8093e-7  &5.9038e-13   &26  &0.1526  &4.9956e-3    &6.2081e-13 \\
~ &16  &5  &0.0557 &6.4339e-11   &5.9274e-13                 &28  &0.1643 &1.1330e-5   &2.9475e-14\\
~ &24  &6  &0.0581 &6.7734e-12   &5.9037e-14                 &30  &0.1665 &3.2785e-8    &8.9360e-13\\\hline
\end{tabular}
\end{table}

\begin{figure}[!ht]
\centering
\subfigure[Multiple solutions with $\lambda = 1$]{
\includegraphics[width=6.5cm]{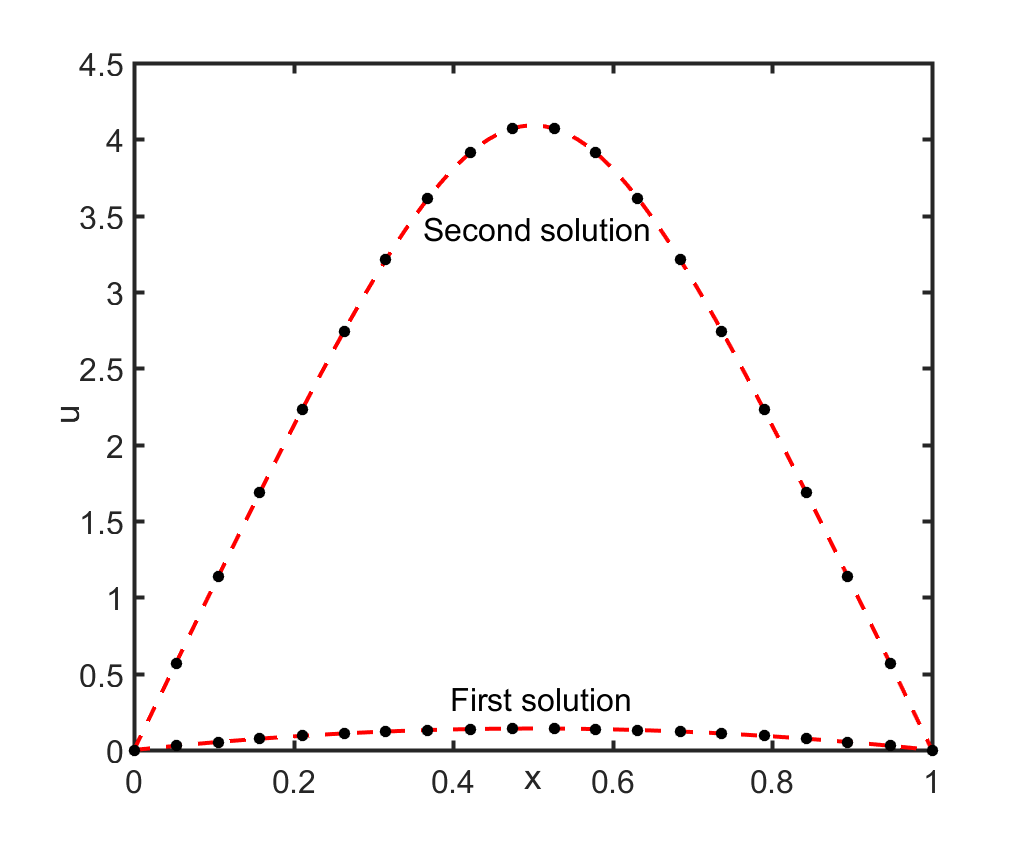}
}
\;
\subfigure[Multiple solutions with $\lambda = 2$]{
\includegraphics[width=6.5cm]{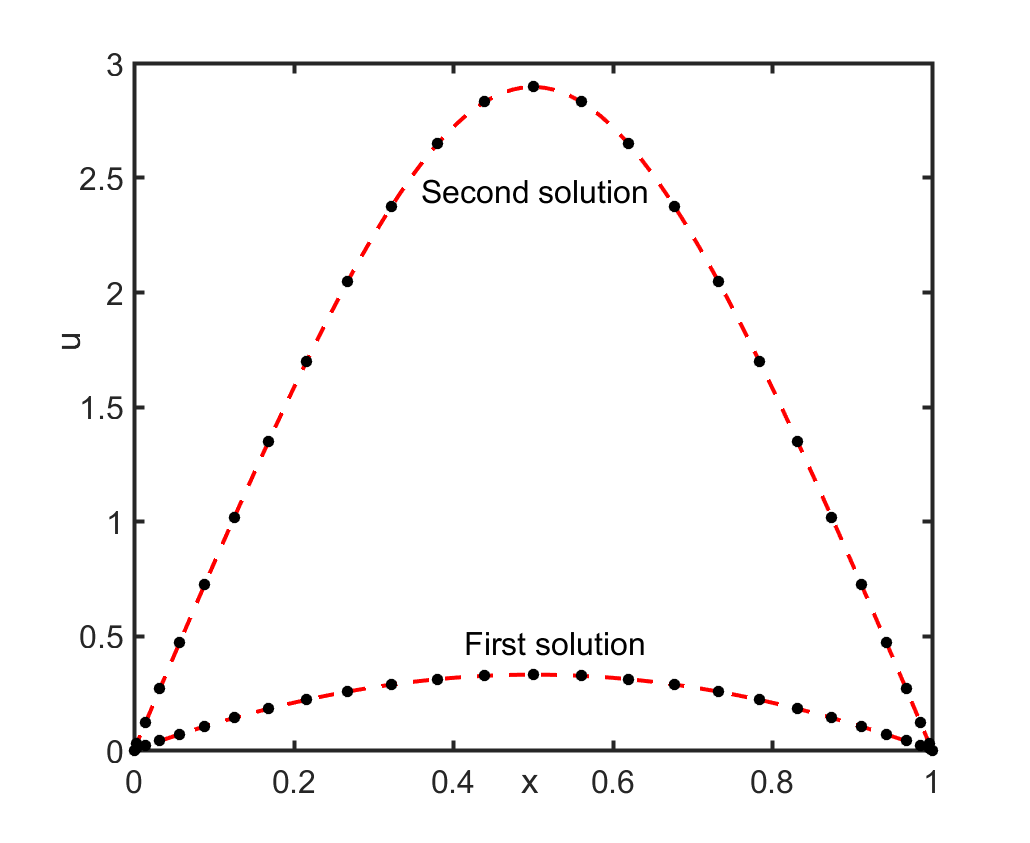}
}
\caption{Multiple solutions of the Bratu-Gelfand model with ${\color{red} \lambda = 1, 2}$.}
\label{Newfigure4.1}
\end{figure}

In Table \ref{NewTable61},  we list the LSTR iteration number, computational time and numerical errors compared with the reference solution (obtained from refined grids) with $\lambda=1,2$ and the initial guess  $-\cos({\rm IG}).$  It indicates clearly that this method finds these two solutions, illustrated in Figure \ref{Newfigure4.1},  very fast  and accurately.

\begin{figure}[!th]
\begin{centering}
\includegraphics[width=6cm,height=5cm]{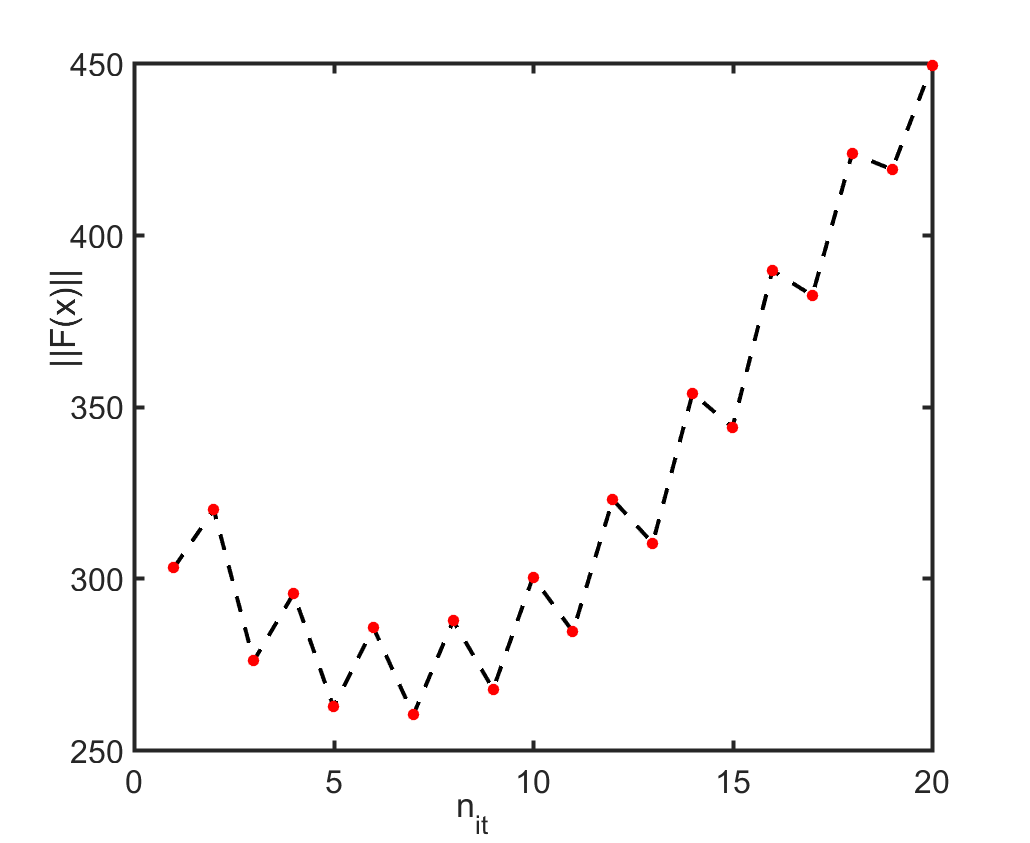}
\caption{The magnitude of $\|{\bs F}({\bs x})\|$ at various iterative steps of the Newton-Deflation method.}\label{NNNFg3.1}
\end{centering}
\end{figure}

In Figure \ref{NNNFg3.1}, we make a similar comparison and record the magnitude of $\|{\bs F}(\bs x^{(i)})\|$ for the first several iterative steps in the Newton-Deflation method.  Again we observe the divergence and sensitivity of the Newton iteration to the initial guess.


\begin{rem}\label{AB} {\em Observe from  {\rm Figure} {\rm \ref{Newfigure4.1}}
that the two solutions satisfy $u'(0)>0$ and $u'(1)<0.$ In fact, it is true for all
 $\lambda\in (0,\lambda^{*}).$ It is evident that by \eqref{3.5} with $F(u)=e^u,$ we have $u''(x)<0,$ so
 $u'(1)<u'(x)<u'(0)$ for all $x\in (0,1).$  Note that if $u'(0)\le 0,$ then $u(x)$ is strictly decreasing in $(0,1),$
 so $u(1)<u(x)<u(0).$  This leads to the contradiction, which implies $u'(0)>0.$ Similarly, we can show  $u'(1)<0.$ \qed
}
\end{rem}

\medskip
\noindent\underline{\bf Case 2:}\; We consider  \eqref{3.5} with  $F(u) = 1 + u^p$ and the boundary conditions:
 $u'(0)=u(1) = 0.$  In the spectral-Galerkin discretization, we use the following basis functions to meet the   boundary conditions
 \begin{equation*}
\psi_k(x) = L_k(x) - \frac{k^2 + 6k + 7}{2(k+2)^2}L_{k+1}(x) - \frac{(k+1)^2}{2(k+2)^2}L_{k+2}(x),
\end{equation*}
which satisfies $\psi'_{k}(-1) = \psi_k(1) = 0$ for any $k$.

 We compare the performance of our proposed method with that of  the bootstrapping
algorithm composed of  the homotopy and two-grid techniques in \cite{hao2014bootstrapping}.
 We quote  \cite[Table 1]{hao2014bootstrapping}  in the left half of  Table \ref{NNewTable61}, where $n_1$ (resp. $n_2)$ is the number of   grid points of the  whole interval (resp. number of grid points in  each of the $n_1-1$ subintervals), together with the numerical errors and computational time for two solutions for $\lambda=1.2$ and $p=4$ in $F(u)$.
 On the right half of Table \ref{NNewTable61}, we list the accuracy and computational time using our method, which shows
 the distinctive advantages of the latter in terms of time-saving  and faster convergence.

\begin{table}[!h]
\centering\small
\caption{\small A comparison between the bootstrapping method and the LSTR-Deflation method with $p = 4$ and $\lambda =1.2$.}
\label{NNewTable61}
\begin{tabular}{|ccccl|cccc|}
	\hline
~ &\multicolumn{4}{c}{Bootstrapping in \cite[Table 1]{hao2014bootstrapping}} &\multicolumn{4}{c|}{The LSTR-Deflation method}\\\cline{2-5}\cline{6-9}
\specialrule{0em}{0pt}{0pt}
$\lambda$ &$(n_1, n_2)$ &I soln err. &II soln err. & Time(s) &$N$  &I soln err. & II soln err. & Time(s)\\\hline   
\multirow{3}{*}{1.2} &(7, 3)  &1.1053e-4  &3.7146e-4 &4m5s  &8  &5.6106e-5  &2.2891e-5 &0.9896s  \\
~ &(21, 2)  &3.2364e-5  &1.1336e-4 &31m16s  &16  &4.4109e-8  &2.4166e-8 &0.3212s  \\
~ &(42, 2)  &8.3840e-6  &3.0290e-5 &5h39m58s  &24  &6.6527e-11  &5.2984e-9 &0.6717s \\\hline
\end{tabular}
\end{table}

\begin{figure}[!ht]
\centering
\subfigure[Multiple solutions with $p = 4$]{
\includegraphics[width=6.5cm]{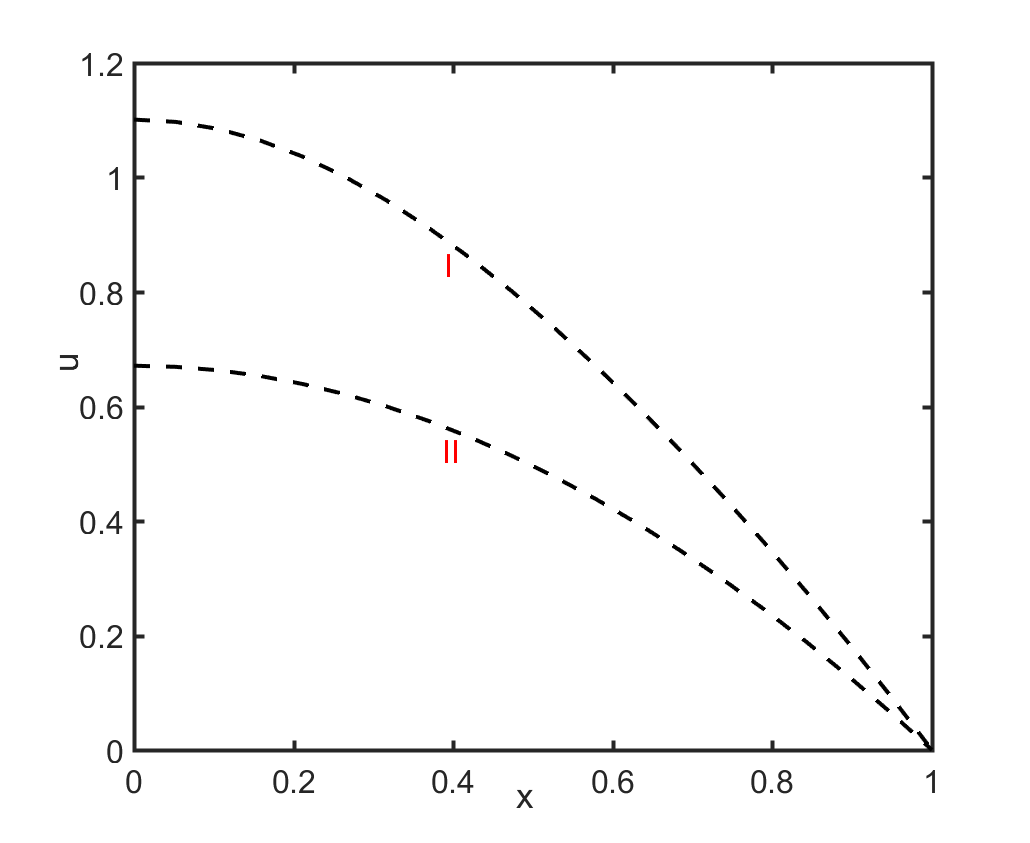}
}
\;
\subfigure[Multiple solutions with $p = 3$]{
\includegraphics[width=6.5cm]{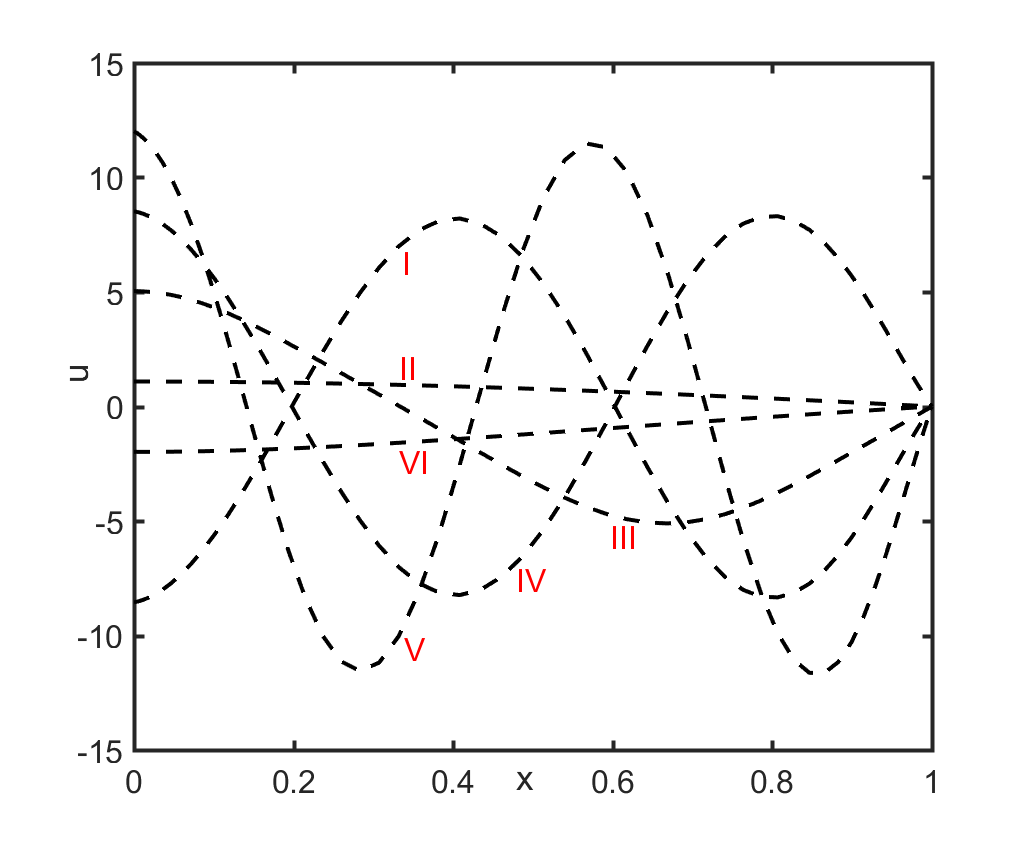}
}
\caption{Multiple solutions for $\lambda = 1.2$.}
\label{Newfigure4.10}
\end{figure}

Notably, the proposed approach allows us to locate more solutions if $p=3$ but $\lambda=1.2$ remains (see Figure \ref{Newfigure4.10} (right)) that seem  not been  reported in literature. We remark that in the deflation process (see \eqref{NNNew2.3.1}),  we have involved  the solutions I, II, and IV to obtain the solution V.

\subsection{PDE Examples}

\subsubsection{Phase separation: Allen-Cahn equation}
As in  Farrell et al. \cite{farrell2015deflation} and \cite{2008Adaptive}, we consider the steady-state Allen-Cahn equation:
\begin{equation}\label{NewExample6}
-\epsilon \Delta u + \epsilon^{-1}(u^{3} - u) = 0  \quad \textrm{in}\;\;\; \Omega = (0, 1)\times(0, 1),
\end{equation}
where the  boundary conditions are as follows
\begin{equation}\label{NNNewExample61}
\begin{dcases}
u = 1 \quad \textrm{on}\;\; \partial\Omega_1:= \{x = 0,1; \; 0 < y < 1\},\\
u = -1  \quad \textrm{on}\;\;  \partial\Omega_2:=\{y = 0,1;\; 0<x<1\}.
\end{dcases}
\end{equation}
Here, $\epsilon$ is
a parameter relating the strength of the free surface tension to the potential term in the free energy.
It is evident that $u=\pm 1$ are solutions to \eqref{NewExample6}, which correspond two different materials.
We can verify that the multiple solutions of \eqref{NewExample6}-\eqref{NNNewExample61} are invariant under a rotation and reflection.
More precisely, if $u(x, y)$ is a solution of \eqref{NewExample6}-\eqref{NNNewExample61}, so is  $-u(1-y, x).$

%
%
%
%

As shown in \eqref{NNNewExample61}, the boundary conditions are not compatible at four corners of $\bar \Omega.$ To facilitate the spectral-Galerkin discretisation,
we properly smooth out the jumps around the corners and render the boundary conditions compatible.
For this purpose, we make a simple linear mapping:  $x\to  2x - 1$, $y\to  2y - 1$ and convert the problem \eqref{NewExample6}-\eqref{NNNewExample61} into $\Omega=(-1,1)^2.$ With a little abuse of notation, we still use $u$ to denote the unknown solution.
Setting $v=u+1$ leads to
\begin{equation}\label{modifyNewExample62}
\begin{dcases}
-4\epsilon \Delta v + \epsilon^{-1}((v - 1)^3 - v + 1)= 0 \quad \textrm{in}\;\; \Omega = (-1, 1)\times(-1, 1),\\
v = 2 \quad \textrm{on}\;\; \partial \Omega_1=\{x = \pm 1,\;  -1<y < 1\},\\
v = 0  \quad \textrm{on}\;\; \partial \Omega_2=\{ y = \pm 1,\; -1<x<1\}.
\end{dcases}
\end{equation}

\begin{figure}[!ht]
\centering
\subfigure[Boundary conditions]{
\includegraphics[width=5.5cm]{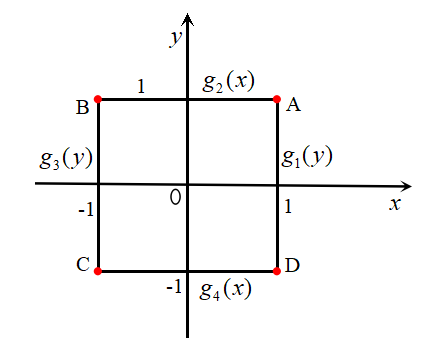}}
\subfigure[$H(y)$ vs $\kappa$]{
\includegraphics[width=5cm]{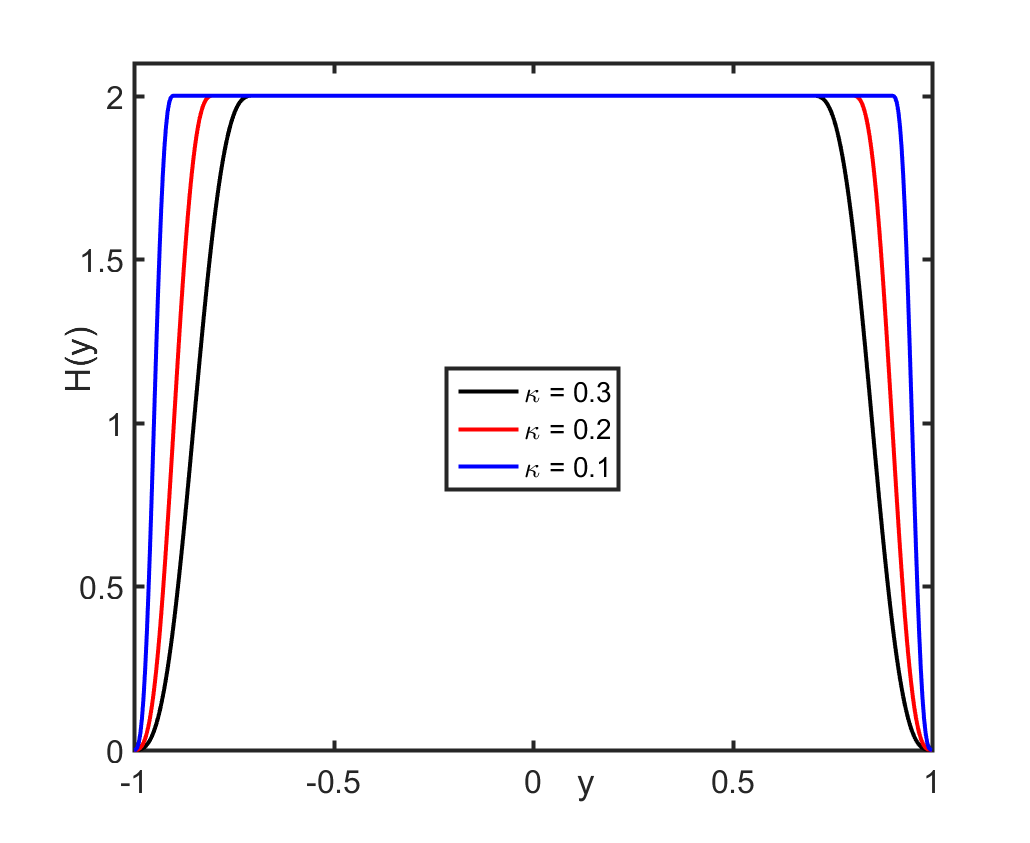}}
\caption{Boundary conditions and $H(y)$ for (\ref{modifyNewExample62}).}\label{vvariblepicture}
\end{figure}
We smooth out $v=2$  by  $H(y)\in C^2[-1,1]$ defined by
\begin{equation}
H(y)=
\begin{dcases}
2\Big[\frac{1}{2} + \frac{1}{2\kappa}(2y + 2 - \kappa) + \frac{1}{2\pi}\sin\Big(\frac{\pi}{\kappa}(2y + 2 - \kappa)\Big)\Big],\;  &
 y \in [-1, -1+\kappa], \\
2, \; & y\in (-1+\kappa, 1-\kappa),  \\
 2 - 2\Big[\frac{1}{2} + \frac{1}{2\kappa}(2y - 2 + \kappa) + \frac{1}{2\pi}\sin\Big(\frac{\pi}{\kappa}(2y - 2 + \kappa)\Big)\Big],\;  & y\in [1-\kappa, 1],
\end{dcases}
\end{equation}
where $\kappa>0$ is a small parameter,  see Figure  \ref{vvariblepicture} (b).   As a result, we derive the following boundary conditions illustrated  in Figure  \ref{vvariblepicture} (a):
\begin{equation}
g_1(y) = g_3(y) = H(y), \quad\quad  g_2(x) = g_4(x) = 0.
\end{equation}
In order to transform the non-homogeneous boundary conditions into homogeneous conditions, we introduce (see \cite{GuoWang2010} for general non-homogeneous boundary conditions)
\begin{equation}\label{nonhom}
\begin{split}
G(x, y) &= \frac{1}{4}\Big\{2(1-y)g_4(x) + 2(1+y)g_2(x) +2(1-x)g_3(y) +2(1+x)g_1(y)\\
 &\quad -(1-x)(1-y)g_3(-1) - (1-x)(1+y)g_2(-1) \\ &\quad - (1+x)(1-y)g_4(1) - (1+x)(1+y)g_1(1)\Big\},
\end{split}
\end{equation}
which satisfies
$$
G(1,y)=g_1(y),\quad G(x,1)=g_2(x),\quad G(-1,y)=g_3(y),\quad G(x,-1)=g_4(x).
$$
Setting $v = w + G$,  we obtain from \eqref{modifyNewExample62} that
\begin{equation}\label{modify61}
\begin{dcases}
-4 \epsilon \Delta w + \epsilon^{-1}[(w+G-1)^3 - (w + G -1)] = 4\epsilon \Delta G \quad  \textrm{in}\;\; \Omega, \\
w = 0 \quad \textrm{on}\;\;  \partial\Omega.
\end{dcases}
\end{equation}
Then we can follow the Legendre spectral-Galerkin method  in Subsection \ref{subsect21} to discretize \eqref{modify61}. We omit the details.

\begin{figure}[!ht]
\begin{center}
\subfigure[I]{ \includegraphics[width=4.5cm]{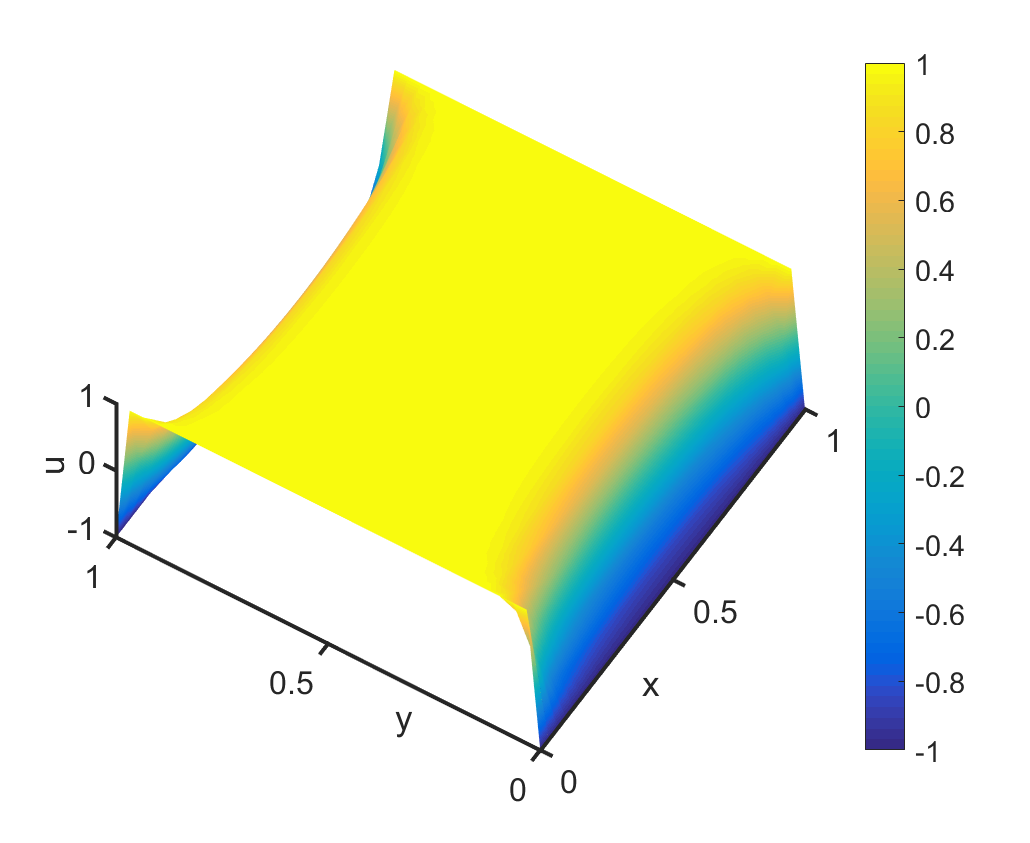}}\;\;
\subfigure[II]{ \includegraphics[width=4.5cm]{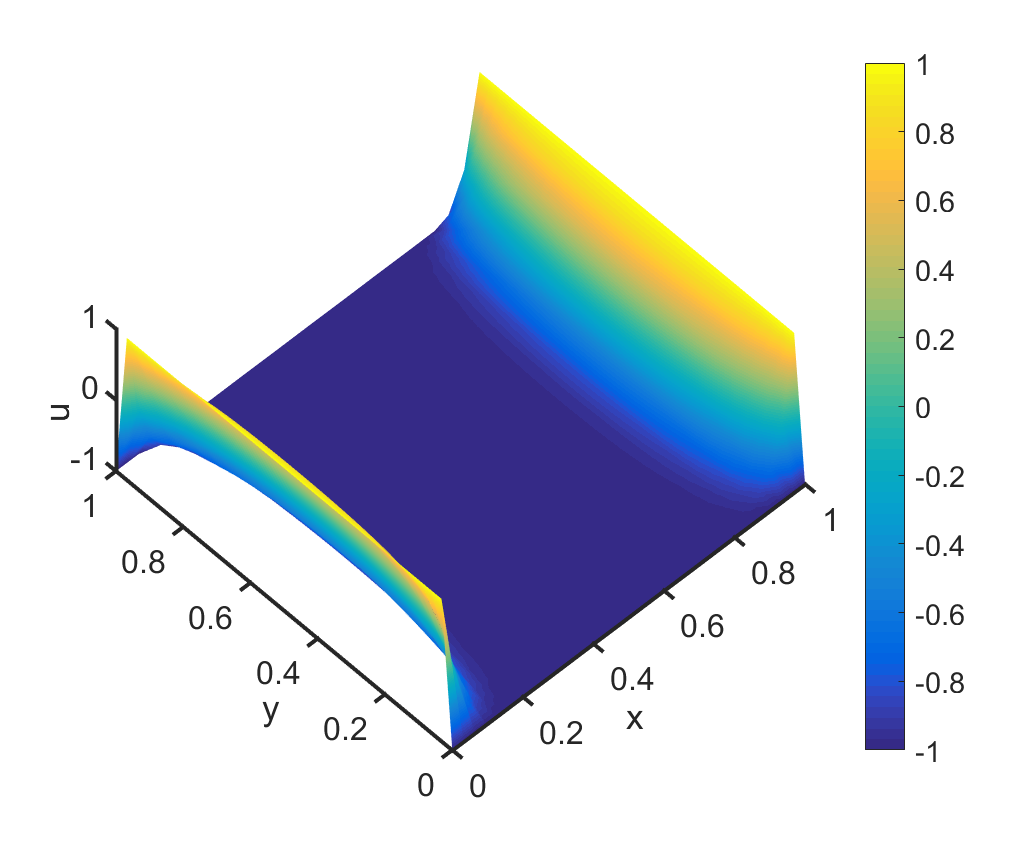}}\;\;
\subfigure[III]{ \includegraphics[width=4.5cm]{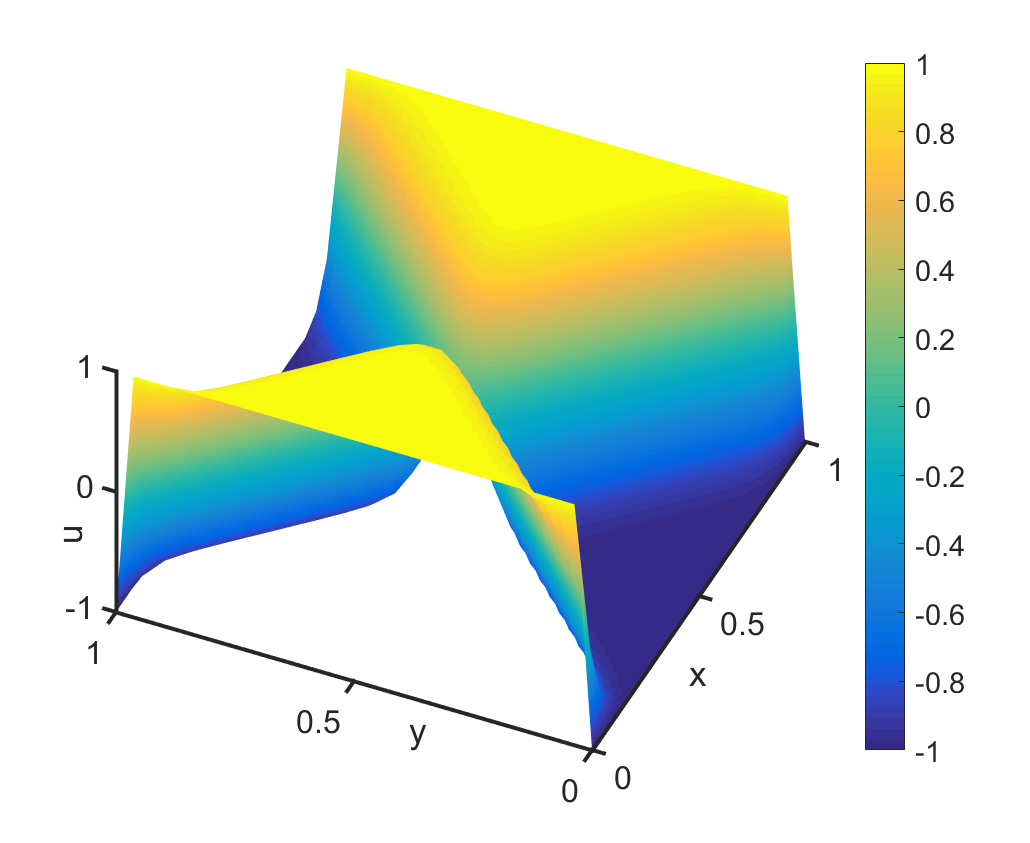}}\\
\caption{Three solutions of the steady Allen-Cahn equation.}
\label{NewExample61}
\end{center}
\end{figure}

In the following tests, we take $\epsilon = 0.04$ and $\kappa = 0.1$. We plot in Figure \ref{NewExample61} the three solutions obtained by our algorithm, whose shapes and patterns are similar to those in \cite{farrell2015deflation, 2018Two}  (however in \cite{2018Two}, only types I and II solutions were reported) without smoothing out the boundary conditions. When $N = 24$, we list the number of iterations, computational time of the LSTR method, the accuracy and the equation errors for each solution in Table \ref{NewTable65}, which shows the high efficiency of the algorithm. It is seen that the first two solutions are identified very fast in a few iterations, but the third solution is relatively expensive and harder to locate.
We also see that the obtained solutions verify the rotational and reflectional property of the solution, where Type-I and Type-II are a pair, and Type-III is invariant under the transforms. Here we have applied three deflations with the same initial guesses to obtain four solutions.

 \begin{table}[!ht]
\centering\small
\caption{\small Performance of  the spectral LSTR-Deflation method.}
\label{NewTable65}\small
\begin{tabular}{ccccc}
	\hline
Solutions &$n_{it}$ &Time(s) &$L^{\infty}$-error &$\|\bs F (\bs x)\|_{\infty}$ \\\hline
\specialrule{0em}{0.5pt}{0.5pt}
I   &39  &5.2710   &3.8109e-10  &1.2984e-13\\\specialrule{0em}{1.5pt}{1.5pt}
II  &23  &4.7952    &6.2910e-11  &3.8652e-13\\\specialrule{0em}{1.5pt}{1.5pt}
III &160  &15.2973    &5.9250e-11  &2.8104e-12\\\specialrule{0em}{1.5pt}{1.5pt}\hline
\end{tabular}
\end{table}

\subsubsection{A model problem analyzed in  Breuer, McKenna and Plum \cite{2003Multiple}}
We consider the following  problem studied in \cite{2003Multiple}:
\begin{equation}\label{NewExample4}
\begin{dcases}
\Delta u + u^{2} = 800\sin(\pi x)\sin(\pi y) \quad&  \textrm{in} \;\; \Omega = (0, 1)\times(0, 1),\\
u = 0  \quad &  \textrm{on} \;\; \partial\Omega.
\end{dcases}
\end{equation}\label{th2}
\begin{figure}[!ht]
\begin{center}
\subfigure[I]{ \includegraphics[width=4cm]{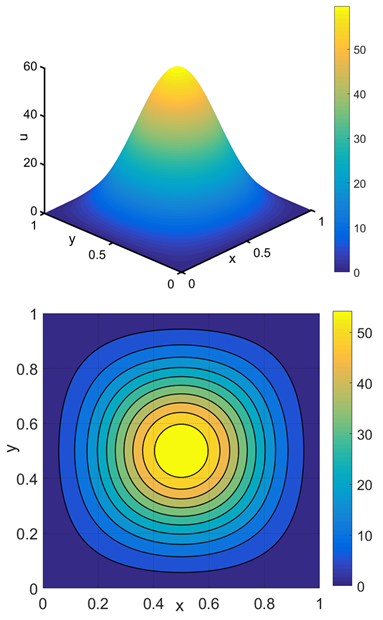}}\;\;
\subfigure[II]{ \includegraphics[width=4cm]{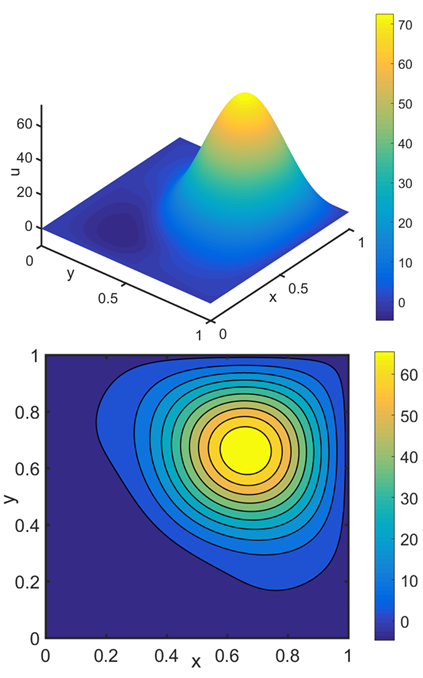}}\;\;
\subfigure[III]{ \includegraphics[width=4cm]{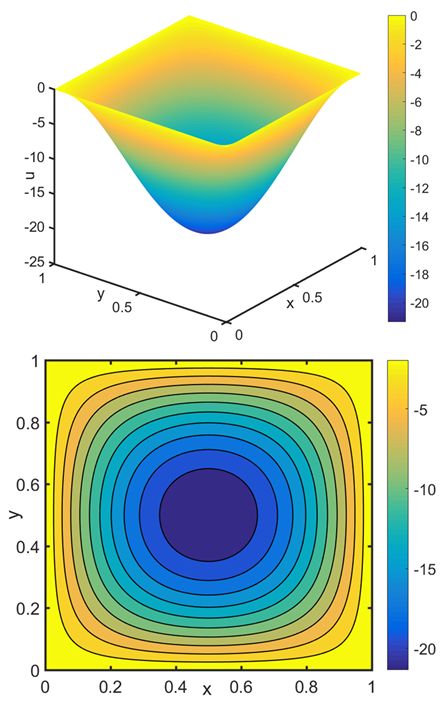}}\;\;
\subfigure[IV]{ \includegraphics[width=6.5cm]{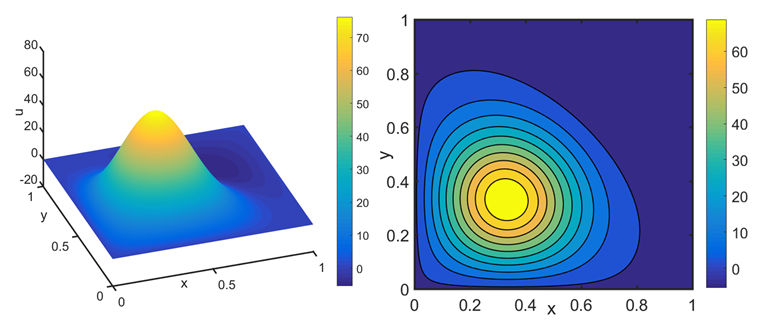}}\;
\subfigure[V]{ \includegraphics[width=6.5cm]{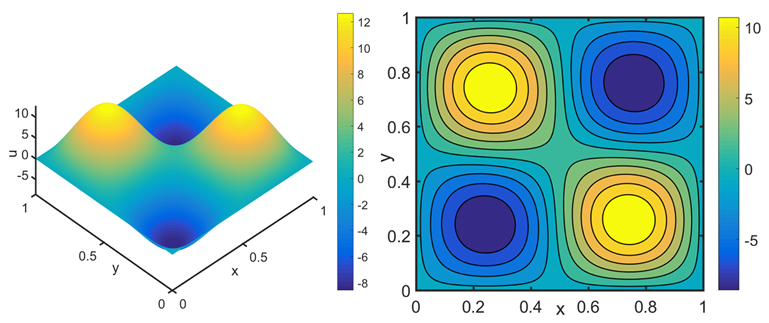}}\\
\caption{Multiple solutions of (\ref{NewExample4}).}
\label{modifyNewExample61}
\end{center}
\end{figure}
 Allgower et al. \cite{allgower2009application} claimed that  (\ref{NewExample4}) has at least four solutions, and these solutions are symmetric (i.e. symmetric with respect to reflections about the axes $x = \frac{1}{2}$, $y = \frac{1}{2}$, $x = y$ and $x = 1 - y$).
 Indeed,  four solutions (i.e. $\textrm{I, II, III}$ and $\textrm{IV}$) obtained from our algorithm  in Figure \ref{modifyNewExample61} satisfy such property, which were
 also presented in \cite{allgower2009application,zhang2013eigenfunction}. More importantly,  we find the fifth solution (i.e.,  $\textrm{V}$  in Figure \ref{modifyNewExample61} (e)), which verifies the claim in \cite{allgower2009application}. Again as shown in Table \ref{NewTable64}, the proposed algorithm has a good performance. Here, it is worth pointing out that when the type-I solution is obtained, the type-II can also be obtained by using the single deflation operator. With the multiple deflation operator in \eqref{NNNew2.3.1}, the type-III solution can be found. While the type-IV and type-V solutions are obtained by using the type-I and type-III solutions with a random perturbation, respectively.

 \begin{table}[!ht]
\centering\small
\caption{\small Numerical results for (\ref{NewExample4}).}
\label{NewTable64}\small
\begin{tabular}{cccccc}
	\hline
Solutions &$n_{it}$ &Time(s) &Symmetry  &$L^{\infty}$-error &$\|\bs F (\bs x)\|_{\infty}$ \\\hline
\specialrule{0em}{0.5pt}{0.5pt}
I   &25  &3.5421  &\makecell{$y = x$, $y = 1 - x$ \\ $x = \frac{1}{2}$, $y = \frac{1}{2}$}  &2.5672e-9  &1.8945e-13\\\specialrule{0em}{1.5pt}{1.5pt}
II  &39  &5.2191  &$y = x$   &6.2145e-7  &5.8230e-14\\\specialrule{0em}{1.5pt}{1.5pt}
III &28  &3.5935  &\makecell{$y = x$, $y = 1 - x$ \\ $x = \frac{1}{2}$, $y = \frac{1}{2}$}    &1.4905e-10  &4.5212e-13\\\specialrule{0em}{1.5pt}{1.5pt}
IV  &30  &4.2418  &$y = x$   &4.6127e-9  &3.9026e-14\\\specialrule{0em}{2.5pt}{3.5pt}
V   &37  &4.6829  &$x = y, y = 1 - x$   &4.6127e-9  &6.8259e-14\\\specialrule{0em}{1.5pt}{1.5pt}\hline
\end{tabular}
\end{table}


\subsubsection{The Henon equation modeling spherical stellar systems} Finally we consider the  Henon equation:
\begin{equation}\label{NNewExample5}
\begin{dcases}
\Delta u + u^{3} = 0 \;\; & \textrm{in} \;\; \Omega=(0,1)^2,\\
u = 0  \quad   & \textrm{on} \;\; \partial\Omega,
\end{dcases}
\end{equation}
as in \cite{chen2004structure}. It is known that  the Henon equation (\ref{NNewExample5}) has infinitely many solutions. Clearly, if $u(x, y)$ is a solution of \eqref{NNewExample5}, so is $-u(x, y)$; furthermore, if rotating a solution $u(x, y)$ by an integer times of 90 degrees around the point ($1/2, 1/2$), then the resulting function, say rot($u(x, y)$), still solves \eqref{NNewExample5}. Thus, in the following tests, we intend to locate some distinct solutions.
\begin{figure}[ht]
  \centering
  \subfigure[I]{\includegraphics[width=4cm]{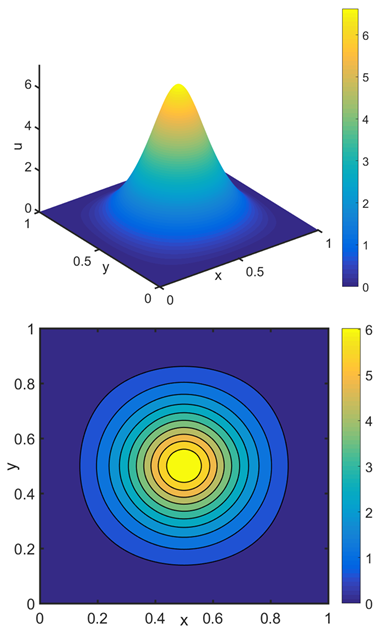}}\;\;
  \subfigure[II]{\includegraphics[width=4cm]{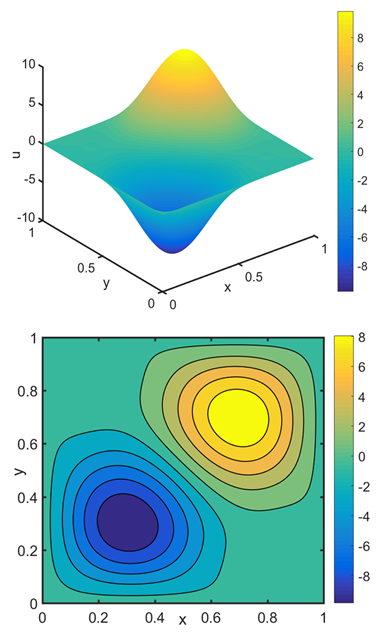}}\;\;
  \subfigure[III]{\includegraphics[width=4cm]{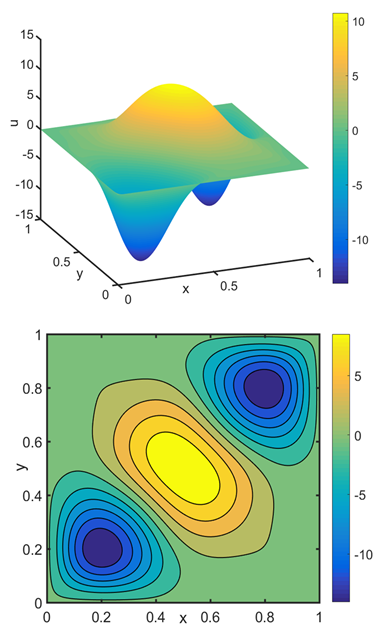}}\\
  \subfigure[IV]{\includegraphics[width=4cm]{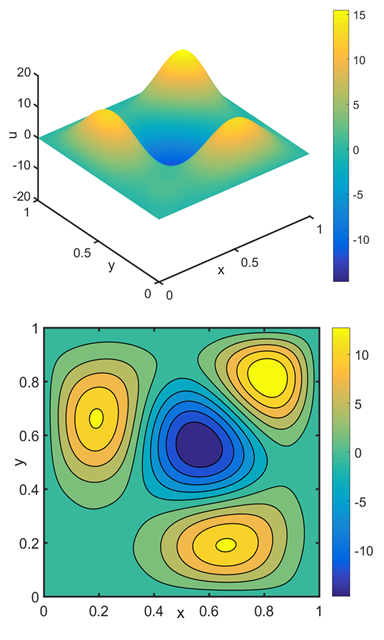}}\;\;
  \subfigure[V]{\includegraphics[width=4cm]{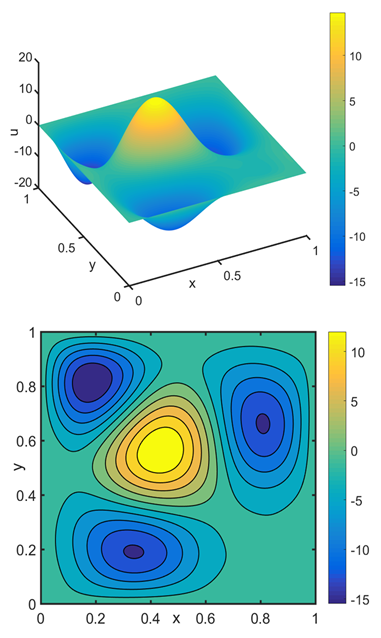}}\;\;
  \subfigure[VI]{\includegraphics[width=4cm]{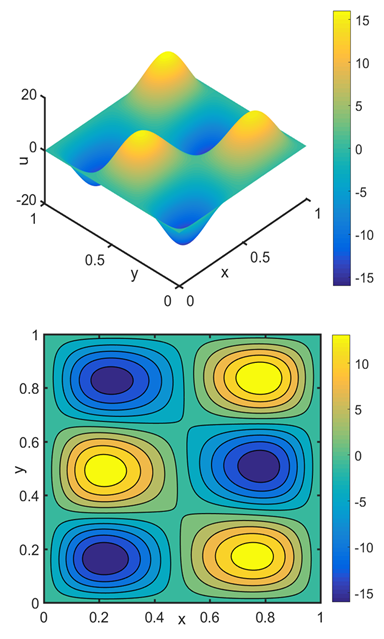}}\;\;\\
  \caption{Six solutions of (\ref{NNewExample5}).}\label{final}
\end{figure}

\begin{table}[!h]
\centering\small
\caption{\small Performance of the spectral LSTR-Deflation methods to \eqref{NNewExample5}.}
\label{NewTable3.2.2.2}
\begin{tabular}{|cccccc|cccccc|}
	\hline
 Types\!\! &$N$  &$n_{it}$ & Time(s)\!\!   &$L^{\infty}$-error &$\|{\bs F}({\bs x})\|_{\infty}$ & Types\!\! &$N$ &$n_{it}$ & Time(s)\!\!  &$L^{\infty}$-error &$\|{\bs F}({\bs x})\|_{\infty}$ \\\hline
\multirow{3}{*}{I} &8  &10  &1.4532 &2.6082e-5   &4.2345e-13 &\multirow{3}{*}{II}     &8  &13  &1.5601 &2.3497e-5 &4.7905e-14\\
~ &16  &15  &1.8902  &3.4508e-7  &9.2409e-13 &~  &16  &15  &1.6012 &5.2109e-7  &6.1208e-14 \\
 ~ &24  &16  &2.0102 &4.8921e-9  &4.7812e-14 &~  &24  &17  &1.8902 &6.0934e-9  &9.5630e-14 \\\hline
 \multirow{3}{*}{III} &8  &20  &1.6703 &2.6071e-4   &7.2409e-13 &\multirow{3}{*}{IV} &8  &21  &1.5608 &4.9021e-5 &9.2301e-13\\
~ &16  &24  &2.0931  &3.4088e-6   &2.0971e-14 &~  &16  &28  &2.0175 &4.0123e-8  &4.8704e-13\\
 ~ &24  &31  &2.3018 &2.0818e-8  &3.0965e-14 &~   &24  &39  &2.3071 &6.0155e-10 &5.2307e-14\\\hline
 \multirow{3}{*}{V} &8  &22  &1.7021 &2.6021e-6   &3.0921e-13 &\multirow{3}{*}{VI}     &8  &24  &1.8902 &5.0905e-4 &8.4302e-14\\
~ &16  &30  &2.0367 &5.8018e-8   &5.7821e-13 &~  &16  &34  &2.3014 &4.9023e-7  &6.8923e-13\\
 ~ &24  &35  &2.9021 &6.2108e-9  &4.0927e-14 &~  &24  &37  &2.7906 &6.6734e-9 &4.0908e-14\\\hline
\end{tabular}
\end{table}

\begin{figure}[ht]
\begin{center}
  \subfigure[VII]{ \includegraphics[width=6.5cm]{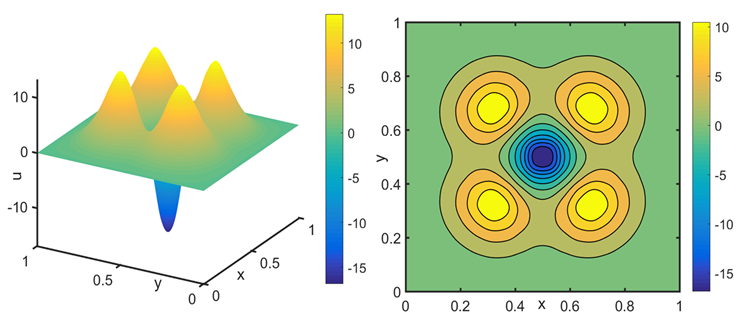}}\;\;
 \subfigure[VIII]{ \includegraphics[width=6.5cm]{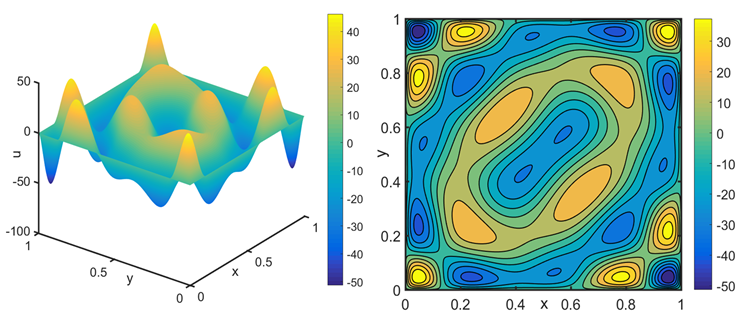}} \\
  \subfigure[IX]{ \includegraphics[width=6.5cm]{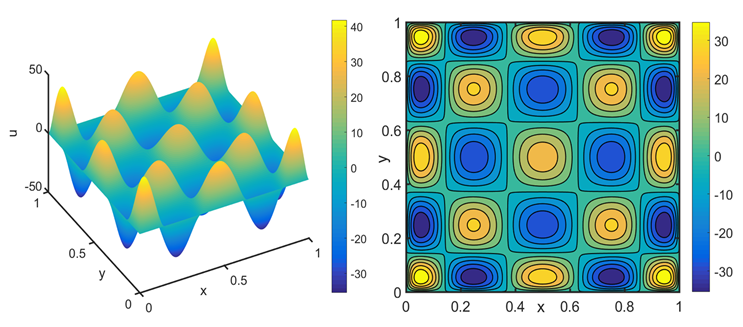}}\;\;
  \subfigure[X]{ \includegraphics[width=6.5cm]{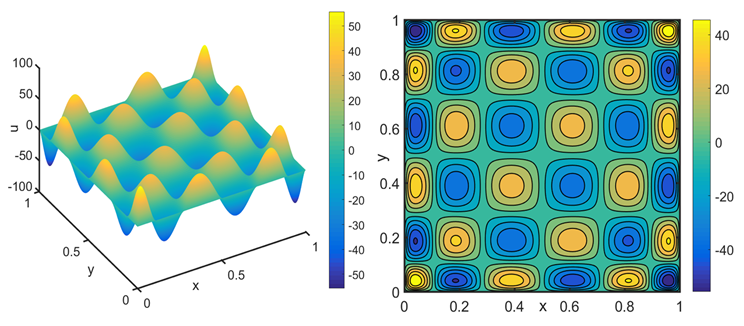}}\\
 \caption{Profiles of multi-peak solutions of (\ref{NNewExample5}).}
   \label{NewExample5}
\end{center}
\end{figure}

Using our  method, we can actually obtain all the multiple solutions presented in \cite{chen2004structure}, but we shall not present them.
Here we intend to demonstrate some solutions which are reported in literature. In Figure \ref{final}, we depict six solutions and their contour plots and in Figure \ref{NewExample5}, we include solutions with   multiple peaks.
It is noteworthy that many more solutions can be obtained from the symmetric and rotational properties. In Tables \ref{NewTable3.2.2.2}-\ref{NewTable3.2.2.1}, we present the outcomes of our algorithm for these solutions with $N = 8, 16$ and 24, where we list iteration number of the LSTR method, the computational time, the maximal norms of  the numerical solutions and the equation errors.

\begin{table}[!h]
\centering\small
\caption{\small Performance of the spectral LSTR-Deflation methods for types VII-X solutions.}
\label{NewTable3.2.2.1}
\begin{tabular}{|cccccc|cccccc|}
	\hline
 Types &$N$  &$n_{it}$ & Time(s)   &$L^{\infty}$-error &$\|{\bs F}({\bs x})\|_{\infty}$ & Types &$N$ &$n_{it}$ & Time(s)  &$L^{\infty}$-error &$\|{\bs F}({\bs x})\|_{\infty}$ \\\hline
\multirow{3}{*}{VII} &8  &18  &1.9029 &3.9014e-5   &8.2305e-14     &\multirow{3}{*}{VIII} &8  &28  &2.4021 &7.0231e-5 &5.0382e-13\\
~ &16  &25  &2.0972 &5.0973e-7   &4.6820e-13 &~  &16  &35  &2.6081 &6.3021e-7  &6.0238e-14\\
 ~ &24  &27  &2.4015 &5.3128e-9  &6.9361e-14 &~  &24  &39  &3.9010 &7.0218e-9  &7.2401e-14\\\hline
 \multirow{3}{*}{IX} &8  &31  &2.4017 &2.4028e-4   &9.4038e-13 &\multirow{3}{*}{X}  &8  &32  &2.0340 &9.3245e-4 &3.9045e-13\\
~ &16  &39  &2.9054 &4.9072e-6   &8.9302e-13 &~  &16  &39  &2.8045 &5.3028e-6  &5.0834e-14\\
 ~ &24  &40  &3.6051 &9.0234e-8  &4.8376e-14 &~  &24  &42  &3.4201 &6.9038e-9 &6.9037e-13\\\hline
\end{tabular}
\end{table}

\section{Concluding remarks}
In this paper,  we proposed an efficient and accurate Legendre spectral-Galerkin method to discretise some differential equations with multiple solutions that lead to nonlinear algebraic systems, and then we combined the nonlinear least-squares and  trust-region method to solve the nonlinear systems and employ the deflation technique to search for their multiple solutions.  We demonstrated that this integrated algorithm has some advantages over certain existing approaches in accuracy, efficiency and capability of finding new solutions. Unfortunately, we didn't provide a rigorous error analysis. Roughly speaking, as the algorithm suggests, it is  quadratically convergent under the conditions in Theorem \ref{thm26} (also applied to the deflated systems through Theorem \ref{thm22}), and the source of the errors is from the spectral-Galerkin discretisation of the original problem. However, much care is needed to glue them together which we hope to provide in a future work.

\medskip
\bigskip

\noindent{\bf Declarations}

\begin{itemize}
  \item {\bf Availability of data and materials:} The data that support the findings of this study are available from the corresponding author upon reasonable request.
  \item {\bf Authors' contributions:} All authors contributed to the study conception and design. The computations and the first draft  were prepared by the first author. All authors read and approved the final manuscript.
\end{itemize}

%
%


%



%
%
%
%
%
%
%
%
%
%
%
%

\bibliographystyle{siam}
\bibliography{myreference1}

\end{document}